\documentclass[dvips,12pt]{article}
\usepackage{graphics}
\usepackage{amsmath}
\usepackage{amssymb}

\edef\savecatcodeat{\the\catcode`@}
\catcode`\@=11

\def\tb@ifSpecChars#1#2{#1}
\def\tb@ifNoSpecChars#1#2{#2}

\def\tableau{%
  \bgroup
  \@ifstar{\let\Tif\tb@ifNoSpecChars\tb@tableauB}
          {\let\Tif\tb@ifSpecChars\tb@tableauB}}

\def\tb@tableauB{
  \@ifnextchar[{\tb@tableauC}{\tb@tableauC[]}}

\def\tb@tableauC[#1]{\hbox\bgroup%
    \let\\=\cr
    \def\bl{\global\let\tbcellF\tb@cellNF}%
    \def\tf{\global\let\tbcellF\tb@cellH}
    \dimen2=\ht\strutbox \advance\dimen2 by\dp\strutbox%
    \ifx\baselinestretch\undefined\relax%
    \else%
       \dimen0=100sp \dimen0=\baselinestretch\dimen0%
       \dimen2=100\dimen2 \divide\dimen2 by\dimen0%
    \fi%
    \let\tpos\tb@vcenter
    \tb@initYoung
    \tb@options#1\eoo
    \let\arrow\tb@arrow%
    \dimen0=\Tscale\dimen2%
    \dimen1=\dimen0 \advance\dimen1 by \tb@fframe%
    \lineskip=0pt\baselineskip=0pt
    \def\tb@nothing{}%
    \def\endcellno{$\rss\egroup\bss\egroup}
    \def\endcell{\endcellno\kern-\dimen0}
    \def\begincell{\vbox to\dimen0\bgroup\vss\hbox to\dimen0\bgroup\hss$}%
    \let\overlay\tb@overlay%
    \let\fl\tb@fl%
    \let\lss\hss\let\rss\hss\let\tss\vss\let\bss\vss
    \def\mkcell##1{
        \let\tbcellF\tb@cellD
        \def\tb@cellarg{##1}
        \ifx\tb@cellarg\tb@nothing\let\tb@cellarg\tb@cellE\fi%
        \begincell\tb@cellarg\endcellno
        \tbcellF}
    \let\savecellF\tbcellF
     \Tif{\catcode`,=4\catcode`|=\active}{}\tb@tableauD}%

\let\tb@savetableauD\tableauD
{
    \catcode`|=\active \catcode`*=\active \catcode`~=\active%
    \catcode`@=\active
\gdef\tableauD#1{%
  \Tif{
    \mathcode`|="8000 \mathcode`*="8000%
    \mathcode`~="8000 \mathcode`@="8000%
    \def@{\bullet}%
    \let|\cr
    \let*\tf
    \let~\sk
  }{}%
  \tpos{\tabskip=0pt\halign{&\mkcell{##}\cr#1\crcr}}%
  \global\let\tbcellF\savecellF
  \egroup
  \egroup}
}
\let\tb@tableauD\tableauD
\let\tableauD\tb@savetableauD
\let\tb@savetableauD\undefined

\def\tb@options#1{\ifx#1\eoo\relax\else\tb@option#1\expandafter\tb@options\fi}

\def\tb@option#1{%
  \if#1t\let\tpos\tb@vtop\fi
  \if#1c\let\tpos\tb@vcenter\fi
  \if#1b\let\tpos\vbox\fi
  \if#1F\tb@initFerrers\fi
  \if#1Y\tb@initYoung\fi
  \if#1s\tb@initSmall\fi
  \if#1m\tb@initMedium\fi
  \if#1l\tb@initLarge\fi
  \if#1p\tb@initPartition\fi
  \if#1a\tb@initArrow\fi
}

\def\tb@vcenter#1{\ifmmode\vcenter{#1}\else$\vcenter{#1}$\fi}

\def\tb@vtop#1{\hbox{\raise\ht\strutbox\hbox{\lower\dimen0\vtop{#1}}}}

\def\tb@initPartition{\def\Tscale{.3}}
\def\tb@initSmall{\def\Tscale{1}}
\def\tb@initMedium{\def\Tscale{2}}
\def\tb@initLarge{\def\Tscale{3}}

\def\tb@initArrow{\dimen2=1.25em}

\def\tb@initYoung{%
  \def\tb@cellE{}
  \let\tb@cellD\tb@cellN
  \def\sk{\global\let\tbcellF\tb@cellNF}}
\def\tb@initFerrers{%
  \def\tb@cellE{\bullet}
  \let\tb@cellD\tb@cellNF
  \def\sk{\bullet}}

\tb@initMedium

\def\tb@sframe#1{%
  \vbox to0pt{
    \vss
    \hbox to0pt{%
      \hss
      \vbox to\dimen1{
        \hrule depth #1 height0pt
        \vss
        \hbox to\dimen1{
          \vrule width #1 height\dimen1
          \hss
          \vrule width #1
          }%
        \vss
        \hrule height #1 depth 0in
        }%
      \kern-\tb@hframe
      }%
    \kern-\tb@hframe}}

\def\tb@hframe{.2pt}\def\tb@fframe{.4pt}\def\tb@bframe{2pt}
\def\tb@cellH{\tb@sframe{\tb@bframe}}       
\def\tb@cellNF{}                            
\def\tb@cellN{\tb@sframe{\tb@fframe}}       
\let\tbcellF\tb@cellN                       

\def\tb@rpad{1pt}
\def\tb@lpad{1pt}
\def\tb@tpad{1.8pt}
\def\tb@bpad{1.8pt}

\def\tb@overlay{\endcell\@ifnextchar[{\tb@overlaya}{\begincell}}
\def\tb@overlaya[#1]{\vbox to\dimen0\bgroup%
  \tb@overlayoptions#1\eoo%
  \tss\hbox to\dimen0\bgroup\lss$}
\def\tb@overlayoptions#1{\ifx#1\eoo\relax\else\tb@overlayoption#1\expandafter\tb@overlayoptions\fi}

\def\tb@overlayoption#1{
  \if#1t\def\tss{\vskip\tb@tpad}\let\bss\vss\fi
  \if#1c\let\tss\vss\let\bss\vss\fi
  \if#1b\def\bss{\vskip\tb@bpad}\let\tss\vss\fi
  \if#1l\def\lss{\hskip\tb@lpad}\let\rss\hss\fi
  \if#1m\let\lss\hss\let\rss\hss\fi
  \if#1r\def\rss{\hskip\tb@rpad}\let\lss\hss\fi
}

\def\tb@fl{\endcell\begincell\vrule depth 0pt width \dimen0 height \dimen0 \endcell\begincell}

\@ifundefined{diagram}{}{
\def\tb@arrowpad{.5}

\newoptcommand{\tb@arrow}{\@ne}[2]{%
  \endcell
   \begingroup%
   \let\dg@getnodesize\tb@getnodesize
   \dg@USERSIZE=#1\relax%
   \ifnum\dg@USERSIZE<\@ne \dg@USERSIZE=\@ne \fi%
   \dg@parse{#2}%
   \dg@label{\tb@draw{#1}{#2}}}

\def\tb@getnodesize#1#2#3#4#5{\dimen3=\tb@arrowpad\dimen2 #4=\dimen3 #5=\dimen3\relax}
\def\tb@getnodesize#1#2#3#4#5{\ifnum#2=0\ifnum#3=0\tb@getnodesizetail{#4}{#5}\else\tb@getnodesizehead{#4}{#5}\fi\else\tb@getnodesizehead{#4}{#5}\fi}
\def\tb@getnodesizetail#1#2{\dimen3=.5\dimen2 #1=\dimen3 #2=\dimen3}
\def\tb@getnodesizehead#1#2{\dimen3=.5\dimen2 #1=\dimen3 #2=\dimen3}

\def\tb@draw#1#2#3#4{%
        \dg@X=0\dg@Y=0\dg@XGRID=1\dg@YGRID=1\unitlength=.001\dimen0%
        \dg@LBLOFF=\dgLABELOFFSET \divide\dg@LBLOFF\unitlength%
        \dg@drawcalc
        \begincell
        \let\lams@arrow\tb@lams@arrow
        \begin{picture}(0,0)\begingroup\dg@draw{#1}{#2}{#3}{#4}\end{picture}%
        \endcell
        \endgroup
        \begincell}
}

\def\tb@lams@arrow#1#2{%
 \lams@firstx\z@\lams@firsty\z@
 \lams@lastx#1\relax\lams@lasty#2\relax
 \lams@center\z@
 \N@false\E@false\H@false\V@false
 \ifdim\lams@lastx>\z@\E@true\fi
 \ifdim\lams@lastx=\z@\V@true\fi
 \ifdim\lams@lasty>\z@\N@true\fi
 \ifdim\lams@lasty=\z@\H@true\fi
 \NESW@false
 \ifN@\ifE@\NESW@true\fi\else\ifE@\else\NESW@true\fi\fi
 \ifH@\else\ifV@\else
  \lams@slope
  \ifnum\lams@tani>\lams@tanii
   \lams@ht\ten@\p@\lams@wd\ten@\p@
   \multiply\lams@wd\lams@tanii\divide\lams@wd\lams@tani
  \else
   \lams@wd\ten@\p@\lams@ht\ten@\p@
   \divide\lams@ht\lams@tanii\multiply\lams@ht\lams@tani
  \fi
 \fi\fi
 \ifH@  \lams@harrow
 \else\ifV@ \lams@varrow
 \else \lams@darrow
 \fi\fi
}

\catcode`\@=\savecatcodeat
\let\savecatcodeat\undefined

\def\part {\vdash}

\def\snake {\rfloor}

\def\la {\lambda}

\def\coeff {{\Big|}}
\def\sab {\vskip .1in}
\def\sabm {\vskip .25in}

\def \sab {\vskip .5in}

\def\ll {{l}}
\def\contained {\subseteq}
\def\V {{\mathcal V}}
\def\H {{\mathcal H}}

\def\wiggle {\tilde}
\def\sz {.35}
\def\one {\text{\sml 1}}
\def\two {\text{\sml 2}}
\def\three {\text{\sml 3}}
\def\four {\text{\sml 4}}
\def\to {{\def\Tscale{\sz} \tableau[Y]{\one}}}
\def\tot {{\def\Tscale{\sz} \tableau[Y]{\one& \two}}}
\def\tto {{\def\Tscale{\sz} \tableau[Y]{\two\\\one}}}    
\def\tabc {{\def\Tscale{\sz} \tableau[Y]{\one& \two&\three}}}
\def\tcab {{\def\Tscale{\sz} \tableau[Y]{\three\\\one&\two}}}    
\def\tbac {{\def\Tscale{\sz} \tableau[Y]{\two\\\one& \three}}}
\def\tcba {{\def\Tscale{\sz} \tableau[Y]{\three\\\two\\\one}}}    

\def\tabcd {{\def\Tscale{\sz} \tableau[Y]{\one& \two&\three & \four}}}
\def\tbacd {{\def\Tscale{\sz} \tableau[Y]{\two\\\one& \three&\four}}}
\def\tcabd {{\def\Tscale{\sz} \tableau[Y]{\three\\\one& \two&\four}}}
\def\tcdab {{\def\Tscale{\sz} \tableau[Y]{\three&\four\\\one& \two}}}
\def\tcbad {{\def\Tscale{\sz} \tableau[Y]{\three\\\two\\\one&\four}}}
\def\tdabc {{\def\Tscale{\sz} \tableau[Y]{\four\\\one&\two&\three}}}
\def\tbdac {{\def\Tscale{\sz} \tableau[Y]{\two&\four\\\one&\three}}}
\def\tdbac {{\def\Tscale{\sz} \tableau[Y]{\four\\\two\\\one&\three}}}
\def\tdabc {{\def\Tscale{\sz} \tableau[Y]{\four\\\one&\two&\three}}}
\def\tdcab {{\def\Tscale{\sz} \tableau[Y]{\four\\\three\\\one&\two}}}
\def\tdcba {{\def\Tscale{\sz} \tableau[Y]{\four\\\three\\\two\\\one}}}

\def\Habc {H^\tabc}
\def\Hcab {H^\tcab}
\def\Hbac {H^\tbac}
\def\Hcba {H^\tcba}
\def\Habcd {{H^\tabcd}}
\def\Hbacd {{H^\tbacd}}
\def\Hcabd {{H^\tcabd}}
\def\Hcdab {{H^\tcdab}}
\def\Hcbad {{H^\tcbad}}
\def\Hdabc {{H^\tdabc}}
\def\Hbdac {{H^\tbdac}}
\def\Hdbac {{H^\tdbac}}
\def\Hdabc {{H^\tdabc}}
\def\Hdcab {{H^\tdcab}}
\def\Hdcba {{H^\tdcba}}

\def\bop {{\vskip .2in \noindent {\bf Proof:} \hskip .2in}}

\def\eop {\hskip .1in $\circ$\vskip .4in}

\newcommand{\Hti}[1]{{\bf H}_{#1}^{-1}}

\newcommand{\AHr}[2]{{\bf H}_{#1}^{#2}}
\newcommand{\AVr}[2]{{\bar{\bf H}}_{#1}^{#2}}
\newcommand{\ontop}[2]{\stackrel{#1}{#2}}
\def\today{\ifcase\month\or
January\or February\or March\or April\or may\or June\or
July\or August\or September\or October\or November\or
December\fi
\space\number\day, \number\year}
\font\sml=cmr6
\numberwithin{equation}{section}
\newtheorem{thm}{Theorem}[section]
\newtheorem{lemma}[thm]{Lemma}
\newtheorem{prop}[thm]{Proposition}

\newtheorem{conjecture}[thm]{Conjecture}
\newtheorem{fact}[thm]{Fact}
\newtheorem{example}[thm]{Example}

\begin{document}

\title{Positivity for special cases of $(q,t)$-Kostka
coefficients and standard tableaux statistics}
\author{Mike Zabrocki\\ Centre de Recherche Math\'ematiques, Universit\'e \\ de
Montr\'eal/LaCIM, Universit\'e de Qu\'ebec \`a Montr\'eal \\ email: {\small\tt
zabrocki@math.ucsd.edu}}
\date{}
\maketitle
\smallskip
{\small MR Subject Number: {05E10}\hfill}

\bigskip

\begin{abstract}
We present two symmetric function operators $H_3^{qt}$ and $H_4^{qt}$
that have the property $H_{3}^{qt} H_{(2^a1^b)}[X;q,t] =
H_{(32^a1^b)}[X;q,t]$
and $H_4^{qt} H_{(2^a1^b)}[X;q,t] =
H_{(42^a1^b)}[X;q,t]$.  These operators are generalizations
of the analogous operator $H_2^{qt}$ and also have expressions
in terms of Hall-Littlewood vertex operators.
We also discuss statistics, $a_{\mu}(T)$ and $b_{\mu}(T)$,
on standard tableaux such that the $q,t$ Kostka polynomials
are given by the sum over standard tableaux of shape $\la$,
$K_{\la\mu}(q,t) = \sum_T  t^{a_{\mu}(T)} q^{b_{\mu}(T)}$
for the case when
when $\mu$ is  two columns or of the form $(32^a1^b)$
or $(42^a1^b)$.
This provides proof of the positivity of the $(q,t)$-Kostka coefficients
in the previously unknown cases of $K_{\la (32^a1^b)}(q,t)$ and 
$K_{\la (42^a1^b)}(q,t)$.  The vertex operator formulas are used to
give formulas for generating functions 
for classes of standard tableaux that generalize the case when
$\mu$ is two columns.
\end{abstract}

\section{Introduction and Notation}

A partition $\la$ is a
weakly decreasing sequence of non-negative integers with
$\la_1 \geq \la_2 \geq \ldots \geq \la_k \geq 0$.  The length $\ll(\la)$
of the partition is the largest $i$ such that $\la_i >0$.  The partition
$\la$ is a partition of $n$ if $\la_1 + \la_2 + \cdots + \la_{\ll(\la)} = n $.
Young diagrams will be drawn using the French notation with the longest row on
the bottom and will be identified with the partition itself by referring to
a partition as a collection of cells.

For every partition $\la$ there is a corresponding conjugate partition denoted
by $\la'$ where  $\la_i' = $ the number of cells in the $i^{th}$ column of
$\la$.  The arm of a cell $s$ in $\mu$ will be denoted by $a_\mu(s)$ and is the
number of cells that lie to the east of $s$ in $\mu$.  The leg, $l_\mu(s)$, is the
number of cells in $\mu$ that are strictly north.

A skew partition is denoted by $\la \slash \mu$, where it is assumed
that $\mu \contained \la$, and represents the
cells that are in $\la$ but are not in $\mu$.
A skew partition $\la \slash \mu$
is said to be a horizontal strip if there is at most
one cell in each column.  Denote the class of horizontal
strips of size $k$ by $\H_k$ so that the notation $\la \slash \mu
\in \H_k$ means that $\la \slash \mu$ is a horizontal strip with $k$
cells.  Similarly, the class of vertical strips (skew partitions
with at most one cell in each row) will be denoted by $\V_k$.

If $\la$ is a partition, then let $\la^r$ denote the partition with
the first row removed, that is $\la^r = (\la_2, \la_3,
\ldots, \la_{\ll(\la)})$.
Let $\la^c$ denote the partition with the first column removed, so that
$\la^c = (\la_1 -1, \la_2 -1, \ldots, \la_{\ll(\la)} -1 )$.
This allows us to define the border of a partition $\mu$ 
to be the skew partition $\mu \slash
\mu^{rc}$.  Define the $k$-snake of a partition $\mu$ to be the
$k$ bottom most right hand cells of the border of $\mu$ (the choice of
the word ``snake'' is supposed to suggest the cells that slink from
the bottom of the partition up along the right hand edge).  We use the
symbol $ht_k(\mu)$ to denote the height of the $k$-snake.  The symbol
$\mu \snake_k$
will be used to represent a partition with the $k$-snake removed
with the understanding that if removing the $k$-snake does not leave
a partition that this symbol is undefined. 

\sabm
\noindent
{\bf Examples}
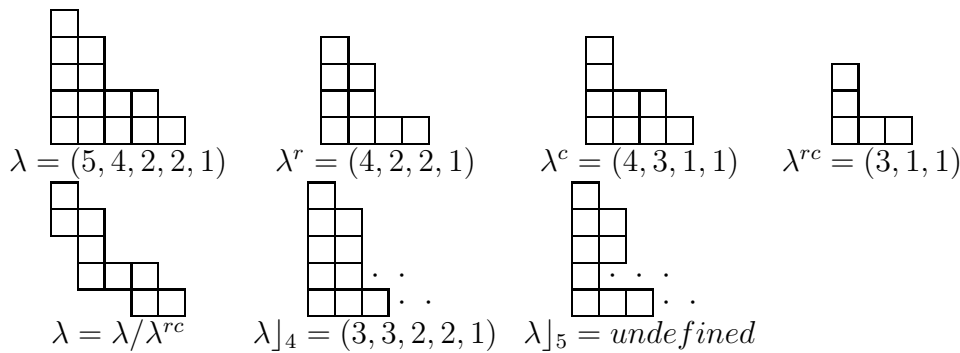
\nopagebreak

\def\Tscale{.7}
\def\ts{{\hskip .1in}}
$$\begin{array}{cccc} 
\ts {{\tableau[Y]{|,|,|,,,|,,,,}}}\ts  &
\ts {{\tableau[Y]{~||,|,|,,,}}}\ts  &
\ts {{\tableau[Y]{~|||,,|,,,}}} \ts &
\ts {{\tableau[Y]{~|~|||,,}}}\ts  \\
\la = (5,4,2,2,1) & \la^r = (4,2,2,1) & \la^c = (4,3,1,1) & \la^{rc} =
(3,1,1)  \\
\ts {{\tableau[Y]{|,|~,|~,,,|~,~,~,,}}}\ts &
\ts {{\tableau[Y]{|,|,|,,~\cdot,~\cdot|,,,~\cdot,~\cdot}}}\ts  &
\ts {{\tableau[Y]{|,|,|,~\cdot,~\cdot,~\cdot|,,,~\cdot,~\cdot}}}\ts & \\
\la = \la \slash \la^{rc} & \la \snake_4 = (3,3,2,2,1) & \la \snake_5 =
undefined &
\end{array}
$$






\sabm
If the shape of $\rho = \la \snake_k$ is given and the height of the
$k$-snake is specified then $\la$ can be recovered
($\la$ is determined from $\rho$ by adding
a $k$-snake of height $h$).  This is because $\la = (\rho_h + k - h +1,
\rho_1+1, \rho_2+1, \ldots, \rho_{h-1}+1, \rho_{h+1}, \rho_{h+2}, \ldots,
\rho_{\ll(\rho)})$.  This will be a partition as long as $k$ is
sufficiently large.

We will consider the ring of symmetric functions in an
infinite number of variables as a subring of
$\mathbb{Q}[x_1, x_2, \ldots]$.  A more precise construction of this
ring can be found in \cite{M} section I.2 (and roughly, the notation
of this reference will be followed).

The Macdonald integral basis \cite{M} for the symmetric functions are defined by
the following two conditions

\begin{align*}
a)& \hskip .1in J_\la = \prod_{s \in \la} (1-q^{a_\la(s)} t^{l_\la(s)+1})
s_\la + \sum_{\mu < \la} s_\mu c_{\mu\la}(q,t) \\
&(s \in \la \text{ \hskip .05in means run over all cells \hskip .05in}  s \text{\hskip .05in in } \la)\\
b)& \hskip .1in \left< J_\la, J_\mu \right>_{qt} = 0 \hskip .1in for \hskip .1in \la \neq \mu 
\end{align*}

\noindent
where $\left< , \right>_{qt}$ denotes the scalar product of symmetric functions defined
on the power symmetric functions by
$\left< p_\la, p_\mu \right>_{qt} = 
\delta_{\la\mu} z_\la \prod_k \frac{1-q^{\la_k}}{1-t^{\la_k}}$
($z_\la$ is the size of the stablizer of the permutations of cycle structure $\la$ and
$\delta_{xy}=1$ if $x=y$ and $0$ otherwise).  The coefficients $c_{\mu\la}(q,t)$
are determined by these two conditions and are rational functions in $q$ and $t$.

  Also define an additional scalar product
$\left< p_\la, p_\mu \right>_{t} = 
\delta_{\la\mu} z_\la \prod_k \frac{1}{1-t^{\la_k}}$.
The $(q,t)$-Kostka coefficients are then given by the expression
$K_{\la\mu}(q,t) = \left< J_\mu[X;q,t], s_{\la}\left[X\right] \right>_{t}$.

We will also refer to the basis
$H_\mu[X;q,t] = \sum_\la K_{\la\mu}(q,t) s_\la[X]$ that is of interest in this paper
as Macdonald symmetric functions.  The $H$ basis is related to the $J$ basis
by a plethystic transformation. Define the basis $H_\mu[X;t] = H_\mu[X;0,t]$ as
Hall-Littlewood symmetric functions.  The $H_\mu[X;t]$ are a transformed version
of the Hall-Littlewood polynomials defined in \cite{M} and are analogous to
the $H_\mu[X;q,t]$.

We will use the notation of $f^\perp$ to denote
the adjoint to multiplication for a symmetric function $f$
with respect to the standard inner product.
Therefore $\left< f^\perp g, h \right> = \left< g, f h \right>$.
Note that $h_k^\perp$ and $e_k^\perp$ act on the Schur function
basis with the formulas
$$e_k^\perp s_\mu = \sum_{\mu\slash\la \in \V_k} s_\la
\hskip .5in h_k^\perp s_\mu = \sum_{\mu\slash\la \in \H_k} s_\la$$


A Schur symmetric function vertex operator is due to Bernstein \cite[p. ~95-6]{M} and
is given by the formula

\begin{equation}
S_m = \sum_{k \geq 0} (-1)^k h_{m+k} e_k^\perp
\end{equation}

It has the property that for $m \geq \mu_1$, $S_m s_\mu[X] = s_{(m,\mu)}[X]$.

There is a Hall-Littlewood symmetric function vertex operator
\begin{equation}
H_m^t = \sum_{i,j\geq 0} t^j (-1)^i h_{m+i+j}[X] e_i^{\perp} h_j^{\perp} =
\sum_{k \geq 0} t^k S_{m+k} h_k^\perp
\end{equation}
with the property that for $m\geq \mu_1$
that $H_m^t H_\mu[X; t] = H_{(m,\mu)}[X;t]$ which is due
to Jing (\cite{J}, \cite{G}, \cite{M}).  This symmetric function operator can
be used to prove the existence of
statistics on column strict tableaux such that $H_\mu[X;t] = \sum_T t^{c(T)}
s_{\la(T)}[X]$.  The action of this operator on the Schur function basis
can be expressed as follows \cite{Za3}.

\begin{prop}\label{snakeruleHm}
Let $\la$ be a partition of $n$, let $m$ be a non-negative integer
and let $k$ be any non-negative
integer such that $m+k \geq \la_1$, then
\begin{equation}
H_m s_\la[X] = \sum_{\mu \slash \la \in \H_{m+k}}
(-1)^{ht_k(\mu)-1} t^{| \la \slash \mu^r |} s_{\mu \snake_k}[X] 
\end{equation}
with the understanding that if $\mu \snake_k$ is not defined then
there is no contribution from that term.
\end{prop}

The object of this research is to find an analogous operator to $H_m^t$
for the Macdonald symmetric functions
$H_\mu[X;q,t]$, $H_m^{qt}$, and use it to derive statistics on standard tableaux
that count the terms in $H_\mu[X;q,t]$.  

Introduce the notation for the operator that acts on symmetric
functions of homogeneous degree $n$ with the formula

\begin{equation}
{\bar H}_m^t = \omega H_m^{\frac{1}{t}} \omega R^t = 
\sum_{i,j\geq 0} t^{n-j} (-1)^i e_{m+i+j}[X] h_i^{\perp} e_j^{\perp}
\end{equation}
where $R^t$ is an operator with the property that $R^t P[X] = t^n P[X]$ for
$P[X]$ a homogeneous polynomial of degree $n$.
In \cite{Za1} we show that

\begin{thm} \label{vop2} The operator
$$H_2^{qt} = H_2^t + q {\bar H}_2^t$$
has the property that
$H_2^{qt} H_{(2^a1^b)}[X;q,t] = H_{(2^{a+1}1^b)}[X;q,t]$.
\end{thm}

In addition we used this operator to show that there exists statistics
on standard tableaux $a_\mu(T)$ and $b_\mu(T)$ for $\mu = (2^a1^b)$
such that
\begin{equation} \label{H2a1btab}
H_{\mu}[X;q,t] = \sum_T q^{b_\mu(T)} t^{a_\mu(T)} s_{\la(T)}[X]
\end{equation}
where the sum is over standard tableaux of size $2a+b$.

This case was already considered by Susanna Fischel in \cite{F} where
it was shown that the $q,t$-Kostka coefficients are a sum over a
subclass of rigged configurations \cite{KR1}, \cite{KR2}.  There exists
a bijection between rigged configurations and standard tableaux, but
the isomorphism between the two sets of objects is not trivial \cite{H},
\cite{KS1}, \cite{KS2}.  This case was also considered in \cite{LM} using
a similar approach but with a different vertex operator for the Macdonald
polynomials.

The problem becomes more difficult when $m=3$ and we will present here
a formula for the vertex operator that adds a row of size $3$, but only works
when it acts on the Macdonald polynomials indexed by less than or equal two columns.

\begin{thm} \label{vop3}
The operator
\begin{align*}
H_3^{qt} &= H_3^t + (e_1[X] H_2^t - H_3^t) q + (e_1[X] {\bar H}_2^t -{\bar H}_3^t) q^2 
+ {\bar H}_3^t q^3\\
&= (1-q)( H_3^t - q^2 {\bar H}_3^t)  + q e_1[X] H_2^{qt}
\end{align*}
has the property that
$H_3^{qt} H_{(2^a1^b)}[X;q,t] = H_{(32^{a}1^b)}[X;q,t]$.
\end{thm}

We will also present an operator that adds a row of size $4$ but again has
the property that it only works when it acts on Macdonald polynomials indexed
by less than or equal two columns.

\begin{thm} \label{vop4}
The operator
\begin{align*}
H_4^{qt} &= H_4^t + (h_1[X] H_3^t - H_4^t) q + (h_2[X] H_2^t - H_4^t) q^2
+ (e_2[X] H_2^t - e_1 H_3^t + H_4^t) q^3 \\
&~~+ (h_2[X] {\bar H}_2^t - h_1 {\bar H}_3^t + {\bar H}_4^t) q^3 
+ (e_2[X] {\bar H}_2^t -{\bar H}_4^t) q^4 + (e_1[X] {\bar H}_3^t -{\bar H}_4^t) q^5
+ {\bar H}_4^t q^6 \\
&=(1-q)(1-q^2)(H_4^t + q^3 \bar{H}_4^t) - q (1-q^2) e_1 (H_3^t - q^2 \bar{H}_3^t) + q^2 (h_2 + q e_2)
(H_2^t + q {\bar H}_2^t)\\
&= (1-q)(1-q^2)(H_4^t + q^3 \bar{H}_4^t) + q (1+q) e_1 H_3^{qt} -
q^2 (e_2 + q h_2) H_2^{qt}
\end{align*}
has the property that
$H_4^{qt} H_{(2^a1^b)}[X;q,t] = H_{(42^a1^b)}[X;q,t]$.
\end{thm}

With the three operators $H_2^{qt}$,
$H_3^{qt}$ and $H_4^{qt}$ acting on $H_{1^b}[X;q,t]$
all of the Macdonald polynomials up to $n=8$ except $(3,3,2)$
are quickly computable
with a computer algebra package such as Maple with John Stembridge's `SF' package or
S\'ebastien Veigneau's package `ACE.'

These two formulas were arrived at by mix of chance, educated guessing, and computer
experimentation.  The proofs of these operators
are not very elegant, but are also not very difficult to follow
and only involve repeated applications of the Macdonald Pieri formula \cite{M}.
They are nice
because they are expressed in terms of the Hall-Littlewood vertex operators
and the action of these operators is well understood \cite{Za3}.  If a formula
for the general Macdonald vertex operator of this sort exists, it will be of the same flavor
as these but it will generalize the Hall-Littlewood vertex rather than have an expression
in terms of it.

In this paper, we will present
statistics on standard tableaux that will show that the Macdonald polynomials 
of the form $H_{(32^a1^b)}[X;q,t]$ and 
$H_{(42^a1^b)}[X;q,t]$, when expanded in terms of
Schur symmetric functions, have coefficients that are polynomials in $q$ and $t$
with non-negative integer coefficients (symmetric functions with this property
will be called Schur positive).

\section{Formulas for $H_{(32^{a}1^b)}[X;q,t]$ and $H_{(42^a1^b)}[X;q,t]$}

An expansion for the Macdonald polynomials $H_{(2^{a}1^b)}[X;q,t]$ in terms
of the Hall-Littlewood symmetric functions was given in \cite{St}.  The coefficients
there are factorable and of a nice form.  Theorem 1.1 in \cite{St}
is the following result

\begin{align}\label{stem}
H_{(2^a1^{b})}[X;q,t] &= \sum_{i=0}^a q^{a-i} (q t^{a+b-i+1}; t)_i
{\frac{(t ; t )_a}{(t ; t)_i(t ; t)_{a-i}}} H_{(2^i1^{b+2a-2i})}[X;t] \\
&= \sum_{i=0}^a c_i^{(a,b)} H_{(2^i1^{b+2a-2i})}[X;t] \nonumber
\end{align}
where
$(a;x)_k = (1-a)(1-ax)\cdots(1-ax^{k-1})$.

Using a translation of the Pieri formula for Macdonald polynomials
\cite[p.~340, eq.~ 6.24(iv)]{M} onto the
$H_\mu[X;q,t]$ basis and that $\mu$ is two columns wide, we can say that
\begin{equation} \label{e1Mac}
e_1[X] H_{(2^{a+1}1^{b})}[X;q,t] = A H_{(32^{a}1^{b})}[X;q,t] + B H_{(2^{a+2}1^{b-1})}[X;q,t]
+ C H_{(2^{a+1}1^{b+1})}[X;q,t]
\end{equation}
where 

\begin{align}
A &= \frac{(1-t^{a+1})(1-q t^{a+b+1}) }{(1 - q^2 t^{a+b+1}) (1-q t^{a+1}) }\\
B &= \frac{(1-t^{b})(1-q ) }{(1 - q t^{b}) (1-q t^{a+1}) }\\
C &= \frac{(1-q)(1-q^2 t^{b}) }{(1 - q t^{b}) (1-q^2 t^{a+b+1}) }
\end{align}

Rearranging this formula to solve for $H_{(32^{a}1^{b})}[X;q,t]$ gives
\begin{equation} \label{H32a1b}
H_{(32^{a}1^{b})}[X;q,t] = A' e_1[X] H_{(2^a1^{b})}[X;q,t] - B' H_{(2^{a+2}1^{b-1})}[X;q,t]
- C' H_{(2^{a+1}1^{b+1})}[X;q,t]
\end{equation}
where 
\begin{align}
A' &= \frac{(1 - q^2 t^{a+b+1}) (1-q t^{a+1}) }{(1-t^{a+1})(1-q t^{a+b+1}) }\\
B' &= \frac{(1 - q^2 t^{a+b+1}) (1-t^{b})(1-q ) }{(1-t^{a+1}) (1 - q t^{b}) (1-q t^{a+b+1}) }\\
C' &= \frac{(1 - q^2 t^{b}) (1-q t^{a+1})(1-q ) }{(1-t^{a+1}) (1 - q t^{b}) (1-q t^{a+b+1}) }
\end{align}

The last two terms on the right hand side of equation (\ref{H32a1b}) can be expanded
in terms of Hall-Littlewood symmetric functions using (\ref{stem}).
The first term on the right hand side of the equation can be expanded
in terms of Hall-Littlewood
symmetric functions by using (\ref{stem}) followed by formula (\ref{e1Mac})
with $q=0$ we have that
\begin{equation} \label{e1HL}
e_1[X] H_{(2^{x}1^{y})}[X;t] =
H_{(2^{x}1^{y+1})}[X;t]
+(1-t^{y}) H_{(2^{x+1}1^{y-1})}[X;t]
+(1-t^{x}) H_{(32^{x-1}1^{y})}[X;t]
\end{equation}

Therefore, arriving at a formula for $H_{(32^{a}1^{b})}[X;q,t]$ is a matter
of applying formulas (\ref{stem}) and (\ref{e1HL}) to (\ref{H32a1b}).  Unfortunately,
the result is not nearly as nice as it was in the case of (\ref{stem}).

\begin{prop} \label{exprH32a1b}
\begin{align}
H_{(32^{a}1^{b})}[X;q,t] &=
\sum_{i=0}^{a} c_{i}^{(a,b)} (1-q^2 t^{a+b+1})  (1-q t^{a+1}) H_{(32^{i}1^{2a+b-2i})}[X;t]
+q^{a+3} H_{(1^{2 a  + b+ 3})}[X;t]  \nonumber \\
&\hskip .2in + \sum_{i=1}^{a+2}  c_{i}^{(a,b)} q
\left(
\frac{(1 - q^2 t^{a+b+1}) (1-q t^{a+1}) }{(1-t^{a+1-{i}}) (1-q t^{a+b+1-{i}}) } \right. \nonumber\\
&\hskip .4in + q 
\frac{(1 - q^2 t^{a+b+1}) (1-q t^{a+1}) (1-t^{i})(1-t^{2 a + b + 4 - 2 i})}
{(1-q t^{a+1+b-{i}}) (1-t^{a+1-{i}}) (1-t^{a+2-{i}} ) (1 - q t^{a+2+b-{i}} )}  \label{bad1} \\
&\hskip .4in - q
\frac{(1 - q^2 t^{a+b+1}) (1-t^{b})(1-q ) (1-t^{a+2})}
{(1-t^{a+2-{i}}) (1-q t^{a+1+b-{i}}) (1-t^{a+1-{i}}) (1 - q t^{b})  } \nonumber \\
&\hskip .4in \left. - \frac{(1 - q^2 t^{b}) (1-q t^{a+1})(1-q ) (1-q t^{a+b+2})}
{(1-q t^{a+b+2-{i}}) (1-q t^{a+1+b-{i}}) (1-t^{a+1-{i}}) (1 - q t^{b}) } \right) \nonumber \\
& \hskip .4in
H_{(2^i1^{2 a + b+3 -2i})}[X;t] \nonumber
\end{align}
\end{prop}

A slightly simpler expression for this result can be given after the vertex
operator is introduced,
but it is not obvious when the coefficients of Proposition \ref{exprH32a1b}
are in a reduced form.

The proof of Proposition \ref{exprH32a1b} uses the following lemmas that each follows
from a simple calculation from the definition of $c_{i}^{(a,b)} = q^{a-i} (q t^{a+b-i+1}; t)_i
{\frac{(t ; t )_a}{(t ; t)_i(t ; t)_{a-i}}}$.

\begin{lemma} For $0 \leq z < x$
\begin{align} 
 \label{coef2}
 c_{z}^{(x,y)} &= c_{z}^{(x,y-1)} \frac{(1-q t^{x+y})}{(1-q t^{x+y-{z}})}\\
 \label{coef3}
 c_{z}^{(x,y)} &= c_{z+1}^{(x,y)} q \frac{(1-t^{z+1})}{(1-t^{x-{z}} ) (1 - q t^{x+y-{z}} )}\\
 \label{coef4}
 c_{z}^{(x,y)} &= c_{z}^{(x-1,y)} q \frac{(1-q t^{x+y}) (1-t^{x})}{(1-q t^{x+y-{z}}) (1-t^{x-{z}})}
\end{align}
\end{lemma}

There are no tricks involved in the reduction of 
(\ref{H32a1b}) to Proposition \ref{exprH32a1b}, just algebraic manipulation.  
Hence we leave the details of the proof to the reader
who may be able to discover a better expression. This proposition is necessary only for
comparison to a similar expression for $H_3^{qt} H_{(2^a1^b)}[X;q,t]$. \eop

To derive a formula for $H_{(42^a1^b)}[X;q,t]$ we will use the same
brute force method for finding equations for the coefficients of
the Hall-Littlewood symmetric functions.  If we add a horizontal
strip of size $2$ on a two column Macdonald summetric function by
multiplying by $g_2 [X;q,0]$ where $g_r[X;q,t]$ is defined
in \cite[eq. 2.8, p. 311]{M} then we have the following
terms

\begin{align}
g_2 [X;q,0] H_{(2^{a+1}1^b)}[X;q,t] =
&A H_{(2^{a+2}1^{b})}[X;q,t] + B H_{(32^{a+1}1^{b-1})}[X;q,t] + 
\nonumber \\ 
&C H_{(32^{a}1^{b+1})}[X;q,t] + D H_{(42^a1^b)}[X;q,t]
\end{align}

\noindent
where $A, B, C,$ and $D$ are given as 

\begin{align}
A &= \frac{1}{(1-q^2 t^{a+b+1}) (1-q t^{a+1})} \\
B &= \frac{(1-t^b)(1-q t^{a+b+1})}{(1-q)(1-q^2 t^{a+1})
(1-q t^b)(1-q^2 t^{a+b+1})} \\
C &= \frac{(1-t^{a+1}) (1-q^2 t^{b})}{(1-q)(1-q t^b)(1-q^3 t^{a+b+1})
(1-q t^{a+1})} \\
D &= \frac{(1-t^{a+1})(1-q t^{a+b+1})}{(1-q)(1-q^2)(1-q^3 t^{a+b+1})
(1-q^2 t^{a+1})}
\end{align}

Rearranging terms, we have that
\begin{align}
H_{(42^a1^b)}[X;q,t] = &A' H_{(2^{a+2}1^{b})}[X;q,t] + 
B' H_{(32^{a+1}1^{b-1})}[X;q,t] +  \nonumber \\
&C' H_{(32^{a}1^{b+1})}[X;q,t] + 
D' g_2 [X;q,0] H_{(2^{a+1}1^b)}[X;q,t]
\label{H42a1bexp}
\end{align}

\noindent
where $A', B', C',$ and $D'$ are given as 

\begin{align}
A' &= -\frac{(1-q)(1-q^2)(1-q^3 t^{a+b+1})(1-q^2 t^{a+1})}
{(1-q^2 t^{a+b+1}) (1-q t^{a+1}) (1-t^{a+1})(1-q t^{a+b+1})} \\
B' &= -\frac{(1-t^b)(1-q^2)(1-q^3 t^{a+b+1})}
{(1-q t^b)(1-q^2 t^{a+b+1})(1-t^{a+1})} \\
C' &= -\frac{ (1-q^2 t^{b})(1-q^2)(1-q^2 t^{a+1})}
{(1-q t^b)(1-q t^{a+1})(1-q t^{a+b+1})} \\
D' &= \frac{(1-q)(1-q^2)(1-q^3 t^{a+b+1})(1-q^2 t^{a+1})}
{(1-t^{a+1})(1-q t^{a+b+1})}
\end{align}

Only the last term in this expression does not have an expression
in terms of Hall-Littlewood symmetric functions yet.
As in the $(32^a1^b)$ case, this can be computed
by translating \cite[p.~340, eq.~ 6.24(i)]{M} and setting $q=0$.
We have that
\begin{align}
h_2[X] H_{(2^x1^{y})}[X;t] = ~
&H_{(2^{x+1}1^y)}+(1-t^{y}) H_{(32^x1^{y-1})}[X;t] + \nonumber \\
&(1-t^{x}) H_{(32^{x-1}1^{y+1})}[X;t] +
(1-t^{x}) H_{(42^{x-1}1^y)}[X;t] \label{h2H2}
\end{align}
and
\begin{align} 
h_1[X]^2 H_{(2^{x}1^{y})}[X;t] =
&H_{(2^{x}1^{y+2})}[X;t] \nonumber\\
&+(2-t^{y+1}-t^{y})H_{(2^{x+1}1^{y})}[X;t]\nonumber\\
&+(1-t^{y})(1-t^{y-1}) H_{(2^{x+2}1^{y-2})}[X;t] \nonumber\\
&+2(1-t^{x}) H_{(32^{x-1}1^{y+1})}[X;t]\nonumber\\
&+(1-t^{y})(2-t^{x+1}-t^{x}) H_{(32^{x}1^{y-1})}[X;t]\nonumber\\
&+(1-t^{x})(1-t^{x-1}) H_{(332^{x-2}1^{y})}[X;t]\nonumber\\
&+(1-t^{x})(1-t) H_{(42^{x-1}1^{y})}[X;t] \label{h1h1H2}
\end{align}
and
\begin{equation}
g_2 [X;q,0] =
\frac{q}{(1-q)(1-q^2)} h_1^2[X] + \frac{1}{1-q^2} h_2[X]
\end{equation}

The next proposition follows by taking coefficients in the formulas above.

\begin{prop} \label{exprH42a1b}
Let $E' = \frac{D'}{(1-q)(1-q^2)}$ and denote the coefficient of
$H_{(32^i1^{2a+b-2i})}[X;t]$ in
$H_{(32^a1^b)}[X;q,t]$ (as given in Proposition \ref{exprH32a1b}) by $d^{(a,b)}_i$
and the coefficient of 
$H_{(2^i1^{3+2a+b-2i})}[X;t]$
in $H_{(32^a1^b)}[X;q,t]$ by $e^{(a,b)}_i$.
The coefficients for the symmetric function $H_{\mu}[X;t]$
in the Macdonald symmetric function $H_{(42^a1^b)}[X;q,t]$
are given by the following table of expressions.

\begin{align}
\mu = (2^i1^{4+2a+b-2i}) & \hskip .4in A' c^{(a+2,b)}_i
+ B' e^{(a+1,b-1)}_i + C' e^{(a,b+1)}_i + 
E'(1-q) c^{(a+1,b)}_{i-1} +  \\
& \hskip .4in
q E' (c^{(a+1,b)}_i
+(2-t^{2a+b-2i+5}-t^{2a+b-2i+4})c^{(a+1,b)}_{i-1}  \nonumber \\
& \hskip 1in +(1-t^{2a+b-2i+6})(1-t^{2a+b-2i+5})c^{(a+1,b)}_{i-2} )  \nonumber\\
\mu = (32^i1^{1+2a+b-2i}) & \hskip .4in B' d^{(a+1,b-1)}_{i}
+ C' d^{(a,b+1)}_{i} +   \\
& \hskip .4inE'(1-q) ((1-t^{2a+b+2-2i}) 
c^{(a+1,b)}_i + (1-t^{i+1}) c^{(a+1,b)}_{i+1}) + \nonumber \\
& \hskip .4inq E' (2(1-t^{i+1})c^{(a+1,b)}_{i+1}
+(1-t^{2a+b+2-2i})(2-t^{i+1}-t^{i}) c^{(a+1,b)}_{i} )
 \nonumber \\
\mu = (332^i1^{2a+b-2i-2}) & \hskip .4in 
q E' (1-t^{i+2})(1-t^{i+1})
c^{(a+1,b)}_{i+2}  \\
\mu = (42^i1^{2a+b-2i}) & \hskip .4in 
E'(1-t^{i+1})(1 -qt) c^{(a+1,b)}_{i+1}      
\end{align}
where if $j<0$ then $c^{(a,b)}_j = 0$ and the coefficients $c^{(a,b)}_i$,
$A'$, $B'$, $C'$ and $D'$ are all given in the text above.  For all $\mu$
that do not follow a pattern in this table, the coefficient of $H_{\mu}[X;t]$
in $H_{(42^a1^b)}[X;q,t]$ is zero.
\end{prop}

Although using this technique for computing expansions of Macdonald polynomials
might be used with other bases, here it seems that a general formula for the
expansion of $H_\mu[X;q,t]$ in terms of $H_\mu[X;t]$ will not be useful since
patterns in the coefficients
do not seem to exist as they did in the case when $\mu = (2^a1^b)$.
We will show that these two particular expansions do have a use to make the
computation of Macdonald polynomials in these special cases much easier.


\section{The Vertex operator $H_3^{qt}$}

Our purpose in giving an expression for $H_{(32^{a}1^{b})}[X;q,t]$ and
$H_{(42^a1^b)}[X;q,t]$, is to
show that the expressions given in Theorems \ref{vop3} and \ref{vop4} have the
properties 
$H_3^{qt} H_{(2^a1^b)}[X;q,t] = H_{(32^{a}1^{b})}[X;q,t]$
and $H_4^{qt} H_{(2^a1^b)}[X;q,t] = H_{(42^a1^b)}[X;q,t]$.

The vertex operators $H_3^{qt}$ and $H_4^{qt}$ are given as expressions using
the $H_m^t$ and ${\bar H}_m^t$ operators and so we will need the action of these
operators on the Hall-Littlewood basis.

\begin{lemma} \label{bH3onHL}
\begin{align*}
{\bar H}_3^t H_{(2^{x}1^{y})}[X;t] &= t^x H_{(2^{x}1^{y+3})}[X;t]
-t^{x+y+1}(1+t) H_{(2^{x+1}1^{y+1})}[X;t] \\
&\hskip .2in -t^{x+y+1}(1-t^y) H_{(2^{x+2}1^{y-1})}[X;t]
+t^{2 x+y+1} H_{(32^{x}1^{y})}[X;t]
\end{align*}
\end{lemma}

\bop
We note the following three commutation relations:
\begin{align}
 {\bar H}_n^t H_m^t &= t^{m-1} H_m^t {\bar H}_n^t \label{comm1}\\
 H_{m-1}^t H_n^t &= t H_m^t H_{n-1}^t + t H_n^t H_{m-1}^t - H_{n-1}^t H_m^t \label{comm3}\\
 H_m^t H_{m+1}^t &= t H_{m+1}^t H_m^t \label{comm2}
\end{align}

A proof of (\ref{comm3}) is in \cite[p.~238]{M} and (\ref{comm2}) is 
a specialization of that.  The proof
of equation (\ref{comm1}) is the same as \cite[Lemma 2.4]{Za1}.

For (\ref{comm3}), in particular we have that
$$H_1^t H_3^t = t H_3^t H_1^t - (1-t) H_2^t H_2^t$$

And ${\bar H}_3 H_{(2^{x}1^{y})}[X;t]$ can be computed by calculating commutation relations.
\begin{align*}
{\bar H}_3 H_{(2^{x}1^{y})}[X;t] &= t^x \left( H_2^t \right)^x {\bar H}_3^t H_{(1^y)}[X;t] \\
&= t^{x} \left( H_2^t \right)^x \left( H_1^t \right)^y {\bar H}_3^t 1 \\
&= t^{x} \left( H_2^t \right)^x \left( H_1^t \right)^y \left( H_{(1^{3})}[X;t] 
-(t+t^2) H_{(21)}[X;t] + t^2 H_{(3)}[X;t] \right) \\
&= t^{x} H_{(2^{x}1^{y+3})}[X;t] - t^{x} (t+t^2) t^{y+1} H_{(2^{x+1}1^{y+1})}[X;t] \\
&\hskip .2in -t^{x+2} t^{y-1} (1-t^y) H_{(2^{x+2}1^{y-1})}[X;t]
 +t^{x+2} t^{x+y-1} H_{(32^{x}1^{y})}[X;t]
\end{align*}
\eop

We now give an outline of the proof of Theorem \ref{vop3}.

\bop
We use equation (\ref{stem}) and Lemma \ref{bH3onHL} 
and show that the expression
is equivalent to Proposition \ref{exprH32a1b}.

\begin{align*}
H_3^{qt} H_{(2^a1^b)}[X;q,t] = &(1-q) H_3^t H_{(2^a1^b)}[X;q,t] + 
q e_1[X] H_{(2^{a+1}1^b)}[X;q,t] \\ 
&+ q^2 (q-1) {\bar H}_3^t H_{(2^a1^b)}[X;q,t] \\
= &(1-q) \sum_{i=0}^{a} c_{i}^{(a,b)} H_{(32^i1^{2a+b-2i})}[X;t] \\
&+ q 
\sum_{i=0}^{a+1} c_{i}^{(a+1,b)} e_1[X] H_{(2^i1^{2 a + 2+b-2i})}[X;t] \\
&+q^2 (q-1) \sum_{i=0}^{a} c_{i}^{(a,b)} \left(
t^i H_{(2^{i}1^{2a+b-2i+3})}[X;t] \right. \\
&\hskip .2in -t^{2a+b-i+1}(1+t) H_{(2^{i+1}1^{2a+b-2i+1})}[X;t]  \\
&\hskip .2in  -t^{2a+b-i+1}(1-t^{2a+b-2i}) H_{(2^{i+2}1^{2a+b-2i-1})}[X;t] \\
&\hskip .2in \left.+t^{2a+b+2} H_{(32^{i}1^{2a+b-2i})}[X;t] ) \right)
\end{align*}

From here it is only algebraic manipulation
to reduce the expression to one like
\begin{align}
=&\sum_{i=0}^{a} c_{i}^{(a,b)} (1-q^2 t^{a+b+1})  (1-q t^{a+1}) H_{(32^{i}1^{2a+b-2i})}[X;t]
+q^{a+3} H_{(1^{2 a  + b+ 3})}[X;t]  \nonumber \\
&+\sum_{i=1}^{a+2} \left(
q c_{i}^{(a+1,b)} 
+ (1-t^{2 a + b + 4 - 2 i}) q c_{i-1}^{(a+1,b)} \right.  \label{besteq} \\
&\hskip .2in + q^2 (q-1) c_{i}^{(a,b)} t^i 
- q^2 (q-1) c_{i-1}^{(a,b)} t^{2a+b+2-i}(1+t) \nonumber \\
&\hskip .2in \left. - q^2 (q-1) c_{i-2}^{(a,b)} t^{2a+b+3-i}(1-t^{2a+b+4-2i}) \right)
 H_{(2^{i}1^{2a+b+3-2i})}[X;t] \nonumber
\end{align}

\noindent
by converting all coefficients in terms of $c_i^{(a,b)}$ using equations (\ref{coef2}),
 (\ref{coef3}) and (\ref{coef4}) we derive the following
expression

\begin{align}=
&\sum_{i=0}^{a} c_{i}^{(a,b)} (1-q^2 t^{a+b+1})  (1-q t^{a+1}) H_{(32^{i}1^{2a+b-2i})}[X;t]
+q^{a+3} H_{(1^{2 a  + b+ 3})}[X;t] \nonumber \\
& +\sum_{i=1}^{a+2} q c_{i}^{(a,b)}  \left(
q \frac{(1-q t^{a+1+b}) (1-t^{a+1})}
  {(1-q t^{a+1+b-{i}}) (1-t^{a+1-{i}})} \right. \nonumber \\
&\hskip .2in + q^2 \frac{(1-t^{2 a + b + 4 - 2 i}) (1-q t^{a+1+b}) (1-t^{a+1}) (1-t^{i})}
  {(1-q t^{a+1+b-{i}}) (1-t^{a+1-{i}}) (1-t^{a+2-{i}} ) (1 - q t^{a+2+b-{i}} )} \label{bad2} \\
&\hskip .2in + q (q-1)  t^i 
- q^2 (q-1)  \frac{(1-t^{i})}{(1-q t^{a+b+1-{i}}) (1-t^{a+1-{i}})}
t^{2a+b+2-i}(1+t) \nonumber \\
& \hskip .2in \left. - q^3 (q-1) \frac{(1-t^{i})(1-t^{i-1}) t^{2a+b+3-i}(1-t^{2a+b+4-2i})}
{(1-q t^{a+b+2-{i}}) (1-t^{a+2-{i}}) (1-q t^{a+b+1-{i}}) (1-t^{a+1-{i}})}
 \right)
 H_{(2^{i}1^{2a+b+3-2i})}[X;t] \nonumber
\end{align}

which can be shown to be equivalent to Proposition \ref{exprH32a1b} by more algebraic
manipulation or by appealing to a computer algebra package such as Mathematica.
\eop

As a corollary, (\ref{besteq}) is an expression for the expansion
of $H_{(32^{a}1^b)}[X;q,t]$ in terms of the basis $H_\mu[X;t]$ since it seems to be
a slightly nicer expression for the Macdonald polynomials than either (\ref{bad1}) or (\ref{bad2}).

We make the following substitution of notation to indicate what the pieces of this operator
represent.  Define the following operators

\begin{equation} \label{defHS3}
\begin{array}{llll}
\Habc &= H_3^t \hskip 1in& \Hbac &= e_1 H_2^t - H_3^t \\
\Hcab &= e_1 \bar{H}_2^t - \bar{H}_3^t \hskip 1in& \Hcba &= \bar{H}_3^t
\end{array}
\end{equation}

Then notice that $H_3^{qt}$ has the expression
\begin{equation} \label{nice3}
H_3^{qt} = \Habc  + q \Hbac + q^2 \Hcab + q^3 \Hcba
\end{equation}

Suddenly the operator that has so far been very unelegant, looks like it is shaping
up.  In the following sections we will show that not only is $H_3^{qt} H_{(2^a1^b)}[X;q,t]$
Schur positive, but so are all of the pieces of this expression
($\Habc H_{(2^a1^b)}[X;q,t]$, $\Hbac H_{(2^a1^b)}[X;q,t]$, $\Hcab H_{(2^a1^b)}[X;q,t]$, and
$\Hcba H_{(2^a1^b)}[X;q,t]$) and that each is a generating function for a subclass of standard
tableaux.

\section{The Vertex operator $H_4^{qt}$}

As with the vertex operator $H_3^{qt}$, we will use what we know about the
action of $H_m^t$ and ${\bar H}_m^t$ on the Hall-Littlewood basis to compute
the action of $(1-q)(1-q^2)(H_4^t + q^3 \bar{H}_4^t) + q (1+q) e_1 H_3^{qt} -
q^2 (e_2 + q h_2) H_2^{qt}$ on equation (\ref{stem}).  The only piece of this
equation that we have not already given an expression for is the action of
$\bar{H}_4^t$ on the Hall-Littlewood symmetric functions.

\begin{lemma}
For $a,b \geq 0$ we have that
\begin{align}
{\bar H}_4^t H_{(2^{a}1^{b})}[X;t] = &t^a H_{(2^{a}1^{b+4})}[X;t] 
- t^{a+b+1} (1+t+t^2) H_{(2^{a+1}1^{b+2})}[X;t] \nonumber\\
&- t^{a+b+1} (1+t-t^{b}-t^{b+1}-t^{b+2}) H_{(2^{a+2}1^{b})}[X;t]  \nonumber\\
&-t^{a+b+1}(1-t^{b-1})(1-t^{b}) H_{(2^{a+3}1^{b-2})}[X;t]  \label{H4bar}\\
&+t^{2 a + b + 2} (1+t) H_{(32^{a}1^{b+1})}[X;t] 
+ t^{2 a + b+ 2} (1+t)(1-t^{b}) H_{(32^{a+1}1^{b-1})}[X;t] \nonumber\\
&+ t^{ 2 a + b + 2 } (1-t^{a}) H_{(332^{a-1}1^b}[X;t]
- t^{2 a + b + 3} H_{(42^{a}1^{b})} \nonumber
\end{align}
\end{lemma}

\bop
First set $a=0$, and prove that for all $b$ greater than or equal to zero we
have that the theorem is true by induction using the identities (\ref{comm1}), (\ref{comm2}),
and (\ref{comm3}) in the special cases of ${\bar H}_4^t H_1^t = t H_1^t {\bar H}_4^t$,
$H_1^t H_4^t = t H_4^t H_1^t + (t^2 - 1) H_3^t H_2^t$,
$H_1^t H_3^t = t H_3^t H_1^t + (t-1) H_2^t H_2^t$, and $H_1^t H_2^t = t H_2^t H_1^t$.
The base case is a calculation of 
$${\bar H}_4^t (1) = H_{(1^4)}[X;t]-t(1+t+t^2) H_{(211)}[X;t]
+t^3 H_{(22)}[X;t] +t^2 (1+t) H_{(31)}[X;t]-t^3 H_{(4)}[X;t]$$

Now for $a>0$ we have that
$${\bar H}_4^t H_{(2^a1^b)}[X;t] = t H_2^t {\bar H}_4^t H_{(2^{a-1}1^b)}[X;t]$$
Using the relations $H_2^t H_3^t = t H_3^t H_2^t$ and $H_2^t H_4^t = H_4^t H_2^t +
(t-1) H_3^t H_3^t$, the result follows.
\eop

We are now ready to prove the vertex operator property.

\bop (of Theorem \ref{vop4})

For each $\mu$ that appears in Proposition \ref{exprH42a1b} we
take the coefficient of $H_\mu[X;t]$ in the expression
\begin{align}
H_4^{qt} H_{(2^a1^b)}[X;q,t] = &(1-q)(1-q^2)(H_4^t + q^3 \bar{H}_4^t) H_{(2^a1^b)}[X;q,t]
+ q (1+q) e_1[X] H_{(32^a1^b)}[X;q,t]  \nonumber\\
& - q^2 (h_1^2[X] + (q - 1) h_2[X]) H_{(2^{a+1}1^b)}[X;q,t]
\end{align}
We need Proposition
\ref{exprH32a1b}, equations
(\ref{stem}), (\ref{e1HL}), (\ref{h2H2}) , (\ref{h1h1H2}) and (\ref{H4bar}), and one
more translation of the Macdonald-Pieri rule.

\begin{align}
e_1[X] H_{(32^{a}1^{b})}[X;t] = &H_{(32^a1^{b+1})}[X;t] +
(1-t^{b}) H_{(32^{a+1}1^{b-1})}[X;t] \nonumber\\
&+(1-t^a) H_{(332^{a-1}1^b)}[X;t]
+(1-t) H_{(42^a1^b)}[X;t]
\end{align}

As in Proposition \ref{exprH42a1b}, denote the coefficient of
$H_{(32^i1^{2a+b-2i})}[X;t]$ in
$H_{(32^a1^b)}[X;q,t]$ (as given in Proposition \ref{exprH32a1b}) by $d^{(a,b)}_i$
and the coefficient of 
$H_{(2^i1^{3+2a+b-2i})}[X;t]$
in $H_{(32^a1^b)}[X;q,t]$ by $e^{(a,b)}_i$.

If $\mu = (42^i1^{2a+b-2i})$, the coefficient is
\begin{equation}
(1-q)(1-q^2)(1 - t^{2 a + b + 3} q^3) c^{(a,b)}_i
+ q (1+q) (1-t) d^{(a,b)}_i+ q^2 (1-t^{i+1})(t-q) c^{(a+1,b)}_{i+i}
\end{equation}
If $\mu = (332^{i}1^{2a+b-2-2i})$, then the coefficient is
\begin{equation}
(1-q)(1-q^2) q^3 t^{ 2a+b+2 } (1-t^{i+1}) c^{(a,b)}_{i+1}
+ q (1+q) (1-t^{i+1}) d^{(a,b)}_{i+1} 
- q^2 (1-t^{i+2}) (1-t^{i+1}) c^{(a+1,b)}_{i+2}
\end{equation}
If $\mu = (32^i1^{2a+b+1-2i})$, then the coefficient is
\begin{align}
&(1-q)(1-q^2) q^3 t^{ 2a+b + 2} (1+t) ( c^{(a,b)}_i 
+  (1-t^{2a+b+2-2i}) c^{(a,b)}_{i-1}) \nonumber\\
&+ q (1+q) ((1-t^{i+1}) e^{(a,b)}_{i+1} +
d^{(a,b)}_{i}+(1-t^{2a+b + 2-2i}) d^{(a,b)}_{i-1}) \\
&- q^2 (2(1-t^{i+1}) c^{(a+1,b)}_{i+1}
+(1-t^{2a+b+2-2i})(2-t^{i+1}-t^{i}) c^{(a+1,b)}_{i} \nonumber\\
&+ (q - 1) ((1-t^{2a+b+2-2i}) c^{(a+1,b)}_{i}
+(1-t^{i+1}) c^{(a+1,b)}_{i+1})) \nonumber
\end{align}
Finally, if $\mu = (2^i1^{2a+b+4-2i})$, then the coefficient is
\begin{align}
&(1-q)(1-q^2) q^3 (t^{i} c^{(a,b)}_{i}
- t^{2a+b+2-i} (1+t+t^2) c^{(a,b)}_{i-1} \nonumber\\
&- t^{2a+b+3-i} (1+t-t^{2a+b+4-2i}-t^{2a+b+5-2i }-t^{2a+b+6-2i }) c^{(a,b)}_{i-2} \nonumber\\
&- t^{2a+b+4-i}(1-t^{2a+b+5-2i})(1-t^{2a+b+6-2i}) c^{(a,b)}_{i-3}) \nonumber\\
&+ q (1+q) (e^{(a,b)}_{i}
+(1-t^{2a+b+5-2i}) e^{(a,b)}_{i-1}) \\
&- q^2 (c^{(a+1,b)}_{i}   
+(2-t^{2a+b+5-2i}-t^{2a+b+4-2i})c^{(a+1,b)}_{i-1}  \nonumber\\
&+(1-t^{2a+b+6-2i})(1-t^{2a+b+5-2i}) c^{(a+1,b)}_{i-2}   \nonumber\\
&+ (q - 1) c^{(a+1,b)}_{i-1}) \nonumber
\end{align}

These can be shown to be equivalent to Proposition \ref{exprH42a1b} by hand (with
an enormous amount of patience) or by reducing these expressions on computer
using a computer algebra package such as Mathematica. \eop

We introduce the following notation as we did in the case of $H_3^{qt}$ to better
demonstrate the structure of this operator.
Define the operators

\begin{equation} \label{defHS4}
\begin{array}{ll}
\Habcd = H_4^t &\hskip .5in
\Hbacd = e_1 H_3^t - H_4^t \\
\Hcdab+\Hcabd = h_2 H_2^t - H_4^t &\hskip .5in
\Hdabc = h_2 {\bar H}_2^t - h_1 {\bar H}_3^t +  {\bar H}_4^t \\
\Hcbad = e_2 H_2^t - e_1 H_3^t +  H_4^t &\hskip .5in
\Hbdac+\Hdbac = e_2 {\bar H}_2^t - {\bar H}_4^t \\
\Hdcab = e_1 {\bar H}_3^t - {\bar H}_4^t & \hskip .5in 
\Hdcba = {\bar H}_4^t
\end{array}
\end{equation}

I am unaware of how to separate the operators
$\Hcdab+\Hcabd$ and $\Hbdac+\Hdbac$, but proofs will work nearly as expected as long
as we consider these as two single entities. Notice that $H_4^{qt}$ has the expression
\begin{align} \label{expv4}
H_4^{qt} = &\Habcd + q \Hbacd + q^2 (\Hcdab + \Hcabd) +
q^3 (\Hdabc + \Hcbad)\\ &~+ q^4 (\Hbdac + \Hdbac) + q^5 \Hdcab + q^6 \Hdcba \nonumber
\end{align}

These operators (and the ones for $H_3^{qt}$) were all defined so that they have the relation
$\omega H^T \omega \coeff_{t \rightarrow 1/t} R^t = H^{\omega T}$ (where $\omega T$ represents
the diagram flipped about its diagonal).

\vskip .2in
Using the relation $H_{\mu'}[X;q,t] = \omega H_\mu[X;t,q]$, we note that
$\omega H_2^{tq} \omega$, $\omega H_3^{tq} \omega$ and $\omega H_4^{tq} \omega$
 add a column of size $2$, $3$ and $4$ respectively to a Macdonald polynomial
indexed by a two row partition.

These vertex operator formulas provide a fast method of computation on
a computer algebra package (such as Maple) of
the Macdonald polynomials for partitions with two rows, two columns, or of the form
$(32^a1^b)$, $(ab1)$, $(42^a1^b)$ or $(ab1^2)$.
Theorem \ref{vop2}, \ref{vop3} and \ref{vop4} are therefore enough to calculate
all of the Macdonald polynomials and the $(q,t)$-Kostka coefficients
through $n=7$ and all but one partition at $n=8$, $(3,3,2)$.

The operators $H_2^{qt}$, $H_3^{qt}$ and $H_4^{qt}$ can be modified so that they work on the other
Macdonald bases.  Let $V$ be the operator with the property
$V J_\mu[X;q,t] = H_\mu[X;q,t]$.  Define $J_m^{qt} = V^{-1} H_m^{qt} V$ then
$J_m^{qt}$ has the corresponding property on the $J_\mu[X;q,t]$ basis.  The formulas
seem to take on an interesting form when expressed in this way.

The methods presented here might be used to handle a few additional cases of deriving formulas
for vertex operators for Macdonald polynomials, but experimental computations indicate that
using the Hall-Littlewood vertex operators will not give nice expressions in general.  There
seems to be another family of vertex operators that has the Hall-Littlewood vertex operators
as a special case that might be used to show
that the $(q,t)$-Kostka coefficients are polynomials with non-negative integer coefficients.

\section{Tableaux, charge, and two column Macdonald polynomials}

A column strict tableau is a diagram of a partition (or skew partition) with
each cell labeled with
a positive integer such that the labels increase weakly
traveling from left to right in the rows
and the labels increase strictly traveling from bottom
to top in the columns.  A standard tableau is a column
strict tableau with the numbers $1$ to  $n$ where $n$ is the
size of the partition.

Let $T$ be a column strict tableau.  Denote the shape
of the tableau by $\la(T)$, the total number of cells in
the diagram by $|T|$, and the
number of cells labeled with an $i$ by the symbol $T_i$.
The content of the tableau will be the tuple $\mu(T) =
(T_1, T_2, \ldots, T_h)$ (where $h$ is the highest label
that appears in the tableau).  T is said to be of
partition content if the content vector $\mu(T)$ is a partition.
The content of a word is defined similarly (the tuple consisting
of (the number of 1's in the word, the number of 2's in the word, etc.))

\sab
We begin by defining an algorithm for calculating the 
statistic called charge on column strict
tableaux of partition weight. This statistic
was introduced by Lascoux and Sch\"utzenberger.

First, charge is defined for words of content weight $\mu=1^n$.
An index is given to each letter in the word.  The index
$0$ is assigned to $1$.  If the letter $i$ has index $k$
then the index of the letter $i+1$ is $k$ if $i+1$ lies to the left of $i$
and the index is $k+1$ if $i+1$ lies to the right of $i$.
The charge of the word is defined to be the sum of the indices.

If $w$ is a word with content of partition weight then it
is first broken up into standard subwords by the following
procedure.
Place an $x$ under the first 1 in the word traveling
from right to left.  Next place an $x$ under the first 2 traveling
to the left from there.  Continue placing an
$x$ underneath each of the letters 1 through $\ll(\mu(w))$ traveling
from right to left and beginning again at the right side of the word
each time the left hand side is reached.
The first standard subword consists of the letters
that have $x$s underneath them read from left to right.
Erase these letters to form a new word $w'$ and repeat the procedure
forming the next standard subword with the labels 1 through $\ll(\mu(w'))$.
Stop when all letters have been erased.  The charge of the word
is then defined to be the sum of the charges of the standard subwords.

The reading word of a tableau is the word formed by reading the
entries in the cells in each of the rows from left to right, 
starting with the top row.  Denote the reading word of $T$ by ${\mathcal R}(T)$.
Lastly, the charge of a tableau T is defined to be the charge
of the reading word of $T$.  Denote the charge of a word $w$ (or tableau $T$) by
$c(w)$ (respectively $c(T)$).

\sabm
\begin{example}
$${{\tableau[sY]{7|3,4,6|2,2,3,5|1,1,1,2,4,8}}}
$$

A tableau of shape $\la = (6,4,3,1)$ and content $\mu = (3,3,2,2,1,1,1,1)$.
The reading word of this tableau is $73462235111248$.  The word has
standard subwords $73625148$, $4231$ and $12$.  The first standard
subword has charge $6$, the second charge $2$, the third charge $1$.
The charge of the tableau is $9$.
\end{example}

If $x$, $y$ and $z$ are letters in the words $u$ and $v$ and 
$w$ and $\wiggle w$ are subwords then say
that two words $u$ and $v$ are elementary Knuth
equivalent if either $u = wxzy{\wiggle w}$ and $v =
wzxy{\wiggle w}$ where $x \leq y < z$ or
$u = wyzx{\wiggle w}$ and $v = wyxz{\wiggle w}$ where
$x < y \leq z$.

Next say that two words $u$ and $v$ are Knuth equivalent
and write $u \sim v$
if they are in the symmetric, transitive, reflexive closure of
the elementary Knuth equivalence.

There are several important facts about words and the charge of words
that we will use to develop statistics here.  These are well known
results and proofs can be found in \cite{LS1}, \cite{B}, \cite{Br},
\cite{Fu}.

\begin{fact} \label{fact0}
Every word is equivalent to the reading word of a unique tableau with partition shape.
\end{fact}

\begin{fact} \label{fact1}
If $u \sim v$ then $c(u) = c(v)$. 
\end{fact}

\begin{fact} \label{fact2}
If $w = w_11^aw_2$ where $w_1$ and $w_2$ are subwords of $w$ that
do not contain $1$ then
$c(w) = c((w_2w_1)\downarrow^1) + |w_2|$ where the notation of
$\downarrow^k$ (dually, $\uparrow_k$)
indicates that the letters of the word to
the left of this symbol have their labels
decreased (increased) by $k$.
\end{fact}

\begin{fact} \label{fact3}
If $w = x_1 x_2 \cdots x_n$ is a standard word, then
$c(w) = \left( \begin{array}{c} n \\ 2\end{array} \right) - c( x_n x_{n-1} 
\cdots x_1)$.
\end{fact}

\begin{fact} \label{fact4}
If $w = x_1 x_2 \cdots x_n$ is a standard word and $w$ is Knuth
equivalent to the reading word of a standard tableau $T$, then
$x_n x_{n-1} \cdots x_1$ is Knuth equivalent to the reading word of
the conjugate tableau $\omega T$.
\end{fact}

The algorithm of Robinson-Schenstead will be used implicitly throughout the following
sections.  Knowledge of the expressions row/column insertion/evacuation and their
relations to Knuth equivalence are assumed in some of the algorithms.  Frequently in
the text we will identify a tableau with its reading word

It is due to the development of this theory
and the results of Lascoux and Sch\"utzenberger \cite{L}, \cite{LS1},
\cite{B}, \cite{Za3}
that the Hall-Littlewood
polynomials $H_\mu[X;t] = \sum_{T \in CST^{\mu}} t^{c(T)} s_{\la(T)}[X]$
where $CST^{\mu}$ is the collection of column strict tableaux of content
$\mu$.

To describe the standard tableaux statisics for the case that the Macdonald
polynomials were indexed by $(2^a1^b)$ we defined
two procedures for building the standard tableaux of size $n+2$ from those of size $n$.
The statement of these definitions will be for arbitrary
$m$ but for now set $m=2$.

Let $T$ be a standard tableaux of shape $\la \vdash n$ and let $\rho$ be a partition of $2n+m$ such
that $\rho \slash \la \in \H_{n+m}$.  Consider the cells that are in $\la$ and not in $\rho^r$ and
perform one column evacuation for each cell from right to left so that cells are evacuated in
increasing order.  Let the row of evacuated cells be $R$ and the remaining tableau be $\tilde T$.
Raise the labels of the cells of $R$ and $\tilde T$ by $m$ and row insert the labels $1$ through $m$
followed by all of $R$ in increasing order.  The result will be the definition of the tableau
$\AHr{m}{\rho} T$. 

We may say that if $\AHr{m}{\rho} T$ and $\rho$ are given, $T$ can be recovered by evacuating 
$12\cdots mR$ a row of size $m+|\la \slash \rho^r|$ from $\AHr{m}{\rho} T$ to leave
$\tilde{T}$ of shape $\rho^r$.  Denote this reverse operation by $(\AHr{m}{\rho})^{-1}$.

\vskip .2in
Also define the operator that adds a column block of $m$ cells to be the transpose
of the operation $\AHr{m}{\rho} T$.  Again, let $T$ be a standard tableaux of shape $\la \vdash n$.
Let $\rho \vdash n+m$ be a partition with $\rho \slash \la \in \V_{n+m}$.  Consider the cells that
are in $\la$ and are not in $\rho^c$ and perform one row evacuation for each cell from top to bottom
so that the cells are evacuated in increasing order.  Let the column of evacuated cells be $C$ and
the remaining tableau be $\tilde T$.  Raise the labels of the cells of $C$ and $\tilde T$ by $m$ and
column insert the labels $1$, through $m$ and all of $C$ in increasing order.  The resulting tableau
will be the definition of $\AVr{m}{\rho} T$.  For the same reason as before,
$T$ can be recovered if $\AVr{m}{\rho} T$ and $\rho$ are given.

\begin{example}
Let $T = \tableau[sY]{2,4|1,3,5,6}$

Let $\rho = (11,3)$ then in the procedure we have that $R = \tableau[sY]{2,4,5}$ and $\tilde{T} =
\tableau[sY]{1,3,6}$.  The labels are all raised by two and $1$, $2$ and $R$ are all inserted
into $\tilde{T}$ to form the tableau

$$\AHr{2}{\rho}T = \tableau[sY]{3,5,8|1,2,4,6,7}$$

Consider $T$ is the same tableau, but now $\rho=(4,3,1^7)$.  In the procedure for adding a
column of size $2$ we have that $C = \tableau[sY]{6}$ and $\tilde{T} = \tableau[sY]{2,4|1,3,5}$.
Therefore,

$$\AVr{2}{\rho} T = \tableau[sY]{8|2,4,6|1,3,5,7}$$
\end{example}

We also introduced the notion of the ${(2^a1^b)}-type$ of
a standard tableau.  The type is defined recursively by `unbuilding' the
standard tableau in blocks of size $2$ (for now set $m=2$ in the following definition).

If $\def\Tscale {.65} \tableau[Y]{1,2, ,m}$ is a subtableau of $T$ then  let $R$ be the first row of
$T$ and $\tilde T$ be $T$ with the first row of $T$ removed.  Define $\Hti{m} T$ to be
the tableau formed by column inserting the labels in $R$ that are not the lables $1$ through $m$
into $\tilde T$ in decreasing order then lower all of the labels in the result by $m$.

If $\def\Tscale {.65} \tableau[Y]{m| |2|1}$ is a subtableau of $T$ then let $C$ be the first column of
$T$ and $\tilde T$ be $T$ with the first column of $T$ removed.  Define $\Hti{m} T$ to be
the tableau formed by row inserting the labels in $C$ that are not a $1$ through $m$
into $\tilde T$ in decreasing order then lower all of the labels in the result by $m$.

Now define the $(2^a1^b)-type$ of $T$ as follows.

\begin{itemize}
\item If $a=0$ and $\mu=(1^b)$ then $type_\mu(T) = (\to^b)$.
\item If $a>1$ and $2$ lies to the right of $1$ in $T$,
then $type_{(2^a1^b)}(T) = (\tot, type_{(2^{a-1}1^b)}(\Hti{2} T))$
\item If $a>1$ and $2$ lies just above $1$ in $T$,
then $type_{(2^a1^b)}(T) = \left(\tto, type_{(2^{a-1}1^b)}(\Hti{2} T) \right)$
\end{itemize}

\begin{example}
Let $T = \tableau[sY]{3|2,6|1,4,5}$.  

$\Hti{2} T = \tableau[sY]{4|2|1,3}$ and $\Hti{2} \Hti{2} T = \tableau[sY]{1,2}$

Therefore, $type_{(2^3)}(T) = \left( \tto, \tto, \tot \right)$.
\end{example}

It was necessary to introduce more machinery in \cite{Za1} to prove the following results.  Here
we will take them as given without defining the handful of tableaux operators that were
necessary to verify that they are true.  The first is Proposition 3.5 and the second is Corollary 3.11
of \cite{Za1}.

\begin{prop}  \label{stats}
Let $n=2a+b$.
The statistic $a_\mu(T)$ and $b_\mu(T)$ where $\mu = (2^a1^b)$ from
equation (\ref{H2a1btab})
may be defined to be 
\begin{align*}
a_\mu(T) &= c(T) - \sum_{i=1}^a ((n+1) - 2i) \chi(type_\mu(T)_i = \tot) \\
b_\mu(T) &= \sum_{i=1}^a \chi \left(type_\mu(T)_i = \tto \right)
\end{align*}
\end{prop}

\begin{prop} \label{type}
Let $n=2a + b$ and $T \in ST^{n}$ and let $\rho$ be a partition of
$2n+2$ such that $\rho \slash \la(T) \in \H_{n+2}$ then
$$type_{(2^{a+1}1^b)}( \AHr{2}{\rho} T) = (\tot, type_{(2^a1^b)}(T))$$
Similarly, if $\rho \slash \la(T) \in \V_{n+2}$ then
$$type_{(2^{a+1}1^b)}( \AVr{2}{\rho} T) = \left(\tto, type_{(2^a1^b)}(T)\right)$$
\end{prop}

\section{Statistics for $\mu = (32^a1^b)$ or $\mu = (42^a1^b)$}

When $H_m^t$ acts on a Macdonald symmetric function, we observe experimentally
that if $m\geq\mu_1-1$, then $H_m^t H_\mu[X;q,t]$ is Schur positive.  This can be explained
for the case that $\mu = (2^a1^b)$ where we have a combinatorial interpretation for the terms.
We then use this result to show that $H^S H_{(2^a1^b)}[X;q,t]$
is Schur positive when $S$ is a standard tableaux of size $3$ or $4$ where the $H^S$ operators
are given in equations (\ref{defHS3}) and (\ref{defHS4}).

The definition of $type_{(2^{a+1}1^b)}( \AHr{2}{\rho} T)$ is 
$(\tot, type_{(2^a1^b)}(\Hti{2} \AHr{2}{\rho} T))$, Proposition \ref{type}
is equivalent to $type_{(2^a1^b)}( \Hti{2} \AHr{2}{\rho} T) = type_{(2^a1^b)}(T)$.
First, extend this result for all $m \geq 2$.

In all of the 
results  in this section
we will assume that $m \geq 2$.

\begin{lemma} \label{mtype}
Let  $T \in ST^{n}$ and let $\rho$ be a partition of $2n+m$ such that
$\rho \slash \la(T) \in \H_{n+m}$ then
$$type_{(2^{a}1^b)}(\Hti{m} \AHr{m}{\rho} T) = type_{(2^{a}1^b)}(T)$$
Similarly, if $\rho \slash \la(T) \in \V_{n+m}$ then
$$type_{(2^{a}1^b)}( \Hti{m} \AVr{m}{\rho} T) = type_{(2^a1^b)}(T)$$
\end{lemma}

\bop
We will just consider the first case since the proof of the second is mostly a matter of changing
$\AHr{k}{\theta}$ to $\AVr{k}{\theta}$ and `row' to `column.'

Let $\tilde{\rho} = \la(\AHr{m}{\rho} T)^r$.  We claim that $\AHr{2}{\tilde{\rho}}
\Hti{m} \AHr{m}{\rho} T = \AHr{2}{ \rho } T$.  

Let $S$ be the tableau of shape $\rho^r$
formed by evacuating the cells in $\la(T)\slash \rho^r$ and let $R$ be the evacuated
cells in the procedures for calculating $\AHr{m}{\rho} T$ and $\AHr{2}{\rho} T$.  We will analyze
the procedures without making changes in the labels for $\AHr{m}{\rho} T$ and instead
of inserting the labels $1$ through $m$, insert $m$ `spaces.' The procedure
for calculating the right hand side  says row insert two spaces and
all of $R$ into $S$.

The procedure for calculating the left hand side says row insert
$m$ spaces and all of $R$ into $S$.
Next, row evacuate the entire first row of the result,
throw away all but $2$ of the spaces and reinsert what is left.
Clearly, since row insertion and row evacuation are inverses of each other,
this is equivalent to the procedure for calculating the right hand side.

The result follows from a couple of applications of Proposition \ref{type}, since
\begin{align}
type_{(2^{a}1^b)}(\Hti{m} \AHr{m}{\rho} T) &=
type_{(2^{a}1^b)}(\Hti{2} \AHr{2}{\tilde{\rho}} \Hti{m} \AHr{m}{\rho} T)  \nonumber \\
&= type_{(2^{a}1^b)}(\Hti{2} \AHr{2}{\rho} T) \\
&= type_{(2^{a}1^b)}( T) \nonumber
\end{align}
\eop

Define an involution, $I_\la^n$, on partitions.
Let $\gamma = \rho \snake_n$. If $\la_h > \gamma_h$ then define
$I_\la^n(\rho) = \gamma$ with an $n-snake$ of
height $h+1$ added. If $\la_h \leq \gamma_h$ then let $I_\la^n(\rho) = \gamma$
with an $n-snake$ of height $h-1$ added.  Lemma 3.15 of \cite{Za1} was the
following result.

\begin{lemma}\label{invol}
 $I_\la^n$ is an involution on partitions $\rho$ such
that $\rho \slash \la \in \H_n$, $\rho \snake_n$ exists and
$\la \neq \rho \snake_n$ with the property that $ht_n(I_\la^n(\rho)) =
ht_n(\rho) \pm 1$ and $\rho \snake_n = I_\la^n(\rho) \snake_n$.
\end{lemma}

\begin{example}
Let $n=10$, $\la=(5,5,2)$, and $\rho=(12,5,5)$, then $\rho \snake_n = (4,4,4)$ and
$ht_n(\rho) = 3$.  $I_\la^n(\rho) = (13,5,4)$ because $I_\la^n(\rho) = (4,4,4)$ with
a $10$ snake of height $2$ added.
\end{example}

\begin{lemma}\label{char}
For a pair $(T, \rho)$ for $T$ a standard tableau, $\rho \slash \la(T) \in \H_{n+m}$
\begin{equation} 
c(\AHr{m}{\rho} T) = c(T) + | \la(T) \slash \rho^r | + \left(\ontop{m}{2}\right)
+ (m-1)n
\end{equation}
\end{lemma}

\bop
We note that $\AHr{m}{\rho} T = (\tilde{T}\uparrow_m)12\cdots m(R\uparrow_m)$ where $\tilde{T}$ and
$R$ are given in the procedure for the calculation of $\AHr{m}{\rho} T$. $T = R\tilde{T}$
and $|R| = |\la(T) \slash \rho^r|$.  Apply Fact \ref{fact2} $m$ times.
\eop

\begin{lemma}\label{tfae}
The following two conditions are equivalent for a pair
$(T,\rho)$ such that $\rho \slash \la(T) \in \H_{n+m}$ and
$\rho \snake_n$ exists.
\begin{align}
\Hti{m} \AHr{m}{\rho} T &= T \\
\la(\AHr{m}{\rho} T) &= \rho \snake_n
\end{align}
Moreover, if $(T,\rho)$ satisfies these these conditions, then $ht_n(\rho)=1$.
\end{lemma}

\bop
Since  $\rho \slash \la(\AHr{m}{\rho} T) \in
\H_{n}$ then if $\la(\AHr{m}{\rho} T) = \rho \snake_n$
then we have that $\rho \slash \rho \snake_n \in \H_n$,
and hence $n=1$.

Let $\tilde{T}$ and $R$ be (respectively) the tableau of shape $\rho^r$ and the row of 
size $|\la(T) \slash \rho^r|$ in the procedure for $\AHr{m}{\rho} T$.

If $\la(\AHr{m}{\rho} T) = \rho \snake_n$ then $\rho \snake_n = (\rho_1-n,\rho^r) = (m+|R|,\rho^r)$
hence the first row of $\AHr{m}{\rho} T$ (which must contain $12\cdots m(R\uparrow_m)$ since it
was row inserted in increasing order) is $12\cdots m(R\uparrow_m)$.  Therefore $\Hti{m} \AHr{m}{\rho} T =
R\tilde{T} = T$.

If $\Hti{m} \AHr{m}{\rho} T = T$ then the first row of $\AHr{m}{\rho} T$ will be
denoted by $12\cdots m\bar{R}$
and the remainder of the tableau will be called $\bar{T}$.  We have that $\bar{R}\bar{T} = T$
and we want to show that $|\bar{R}| = |R|$.  By Lemma \ref{char} and repeated application
of Fact \ref{fact2} we have that

\begin{equation}
c(T)
= c(\bar{T}12\cdots m\bar{R}) - | R | - \left(\ontop{m}{2}\right) - (m-1)n
= c(T) + |\bar{R}| - |R|
\end{equation}
Therefore $\la(\AHr{m}{\rho} T) = (m+|R|,\rho^r) = \rho \snake_n$.
\eop

Call a pair that satisfies the conditions of Lemma \ref{tfae} {\it stable} and a pair
$(T,\rho)$ such that $\rho \slash \la(T) \in \H_{n+m}$ and
$\rho \snake_n$ exists and that does not satisfy the conditions of the lemma, {\it unstable}.
A pair $(T,\rho)$ such that $\rho \slash \la(T) \in \H_{n+m}$ and
$\rho \snake_n$ does not exist will be called {\it immaterial}.

\begin{lemma} \label{Tpinvol}
There exists an involution ${\mathcal I}^{n,m}$ on pairs $(T,\rho)$ 
that are unstable.
The corresponding pair $(\hat{T}, \tilde{\rho}) = {\mathcal I}^{n,m}(T,\rho)$ has the
property that $\AHr{m}{\tilde{\rho}}\hat{T} = \AHr{m}{\rho}T$,
$ht_n(\tilde{\rho}) = ht_n(\rho) \pm 1$ and 
$type_{(2^a1^b)}(T) = type_{(2^a1^b)}(\hat{T})$.
\end{lemma}

\bop
Begin by setting $\la = \la(\AHr{m}{\rho}T)$ and $\tilde{\rho} = I_\la^n(\rho)$ then
the property that $ht_n(\tilde{\rho}) = ht_n(\rho) \pm 1$ follows from Lemma
\ref{invol}.  We want to define $\hat{T} = (\AHr{m}{\tilde{\rho}})^{-1} \AHr{m}{\rho}T$,
but it is not always possible to apply $(\AHr{m}{\tilde{\rho}})^{-1}$ to a tableau.  The
first $m$ columns of $\rho$ and $\tilde{\rho}$ must be the same because $\rho_1$
and $\tilde{\rho}_1$ are both at least $n+m$ and the leftmost $m$ columns of $\rho \snake_n=
\tilde{\rho} \snake_n$
are the same as the leftmost $m$ columns of both $\rho$ and $\tilde{\rho}$.  Therefore the procedure
for applying
$(\AHr{m}{\tilde{\rho}})^{-1}$ to the tableau $\AHr{m}{\rho}T$ evacuates the labels
$12\cdots m$ and hence there is no problem defining $\hat{T}$ in this way.

This is an involution since $I_\la^n(\tilde{\rho})=\rho$ and
$T = (\AHr{m}{\rho})^{-1} \AHr{m}{\tilde{\rho}} \hat{T}$.

We also have that $\AHr{m}{\tilde{\rho}}\hat{T} = \AHr{m}{\rho}T$, and hence by
Lemma \ref{mtype}
\begin{equation}
type_{(2^a1^b)}(T) = type_{(2^a1^b)}( \Hti{m} \AHr{m}{\rho}T)
= type_{(2^a1^b)}( \Hti{m} \AHr{m}{\tilde{\rho}}\hat{T})
= type_{(2^a1^b)}( \hat{T})
\end{equation}

\eop

\begin{example}
Let $m=2$ and $T = \tableau[sY]{5|4|1,2,3}$ and $\rho = (8,3,1)$.

We calculate that
$\AHr{2}{\rho}T = \tableau[sY]{6|3,4|1,2,5,7}$ and find that $(T,\rho)$ is an unstable
pair.  $\tilde{\rho} = (7,4,1)$ and $\hat{T} = (\AHr{2}{\tilde{\rho}})^{-1} \AHr{2}{\rho}T
= \tableau[sY]{4|1,2,3,5}$.

Notice that the types of $T$, $\hat{T}$, and $\Hti{2}\AHr{2}{\rho}T$ are all the same.
$type_{(2,2,1)}(T) = type_{(2,2,1)}(\hat{T}) = type_{(2,2,1)}(\Hti{2}\AHr{2}{\rho}T) = ( \tot, \tot,\to )$.
\end{example}

\vskip .2in
\begin{prop} \label{Hmstats}
\begin{equation*}
H_m^t H_{(2^a1^b)}[X;q,t] = \sum_T q^{b_{(m2^a1^b)}(T)}
t^{a_{(m2^a1^b)}(T)} s_{\la(T)}[X]
\end{equation*}
where the sum is over standard tableau $T$ of size $2a+b+m$ that 
contain $\def\Tscale{.65} \tableau[Y]{1,2, ,m}$ as
a subtableau and $ b_{(m2^a1^b)}(T) = b_{(2^a1^b)}(\Hti{m}T)$ and
\begin{align*}
a_{(m2^a1^b)}(T) = &c(T) - \left( \ontop{m}{2} \right) - (m-1)(2a+b) \\
&- \sum_{i=1}^a ((2a+b+1) - 2i) \chi(type_{(2^a1^b)}(\Hti{m}T)_i = \tot)
\end{align*}
\end{prop}

\noindent
{\bf Remark:} For the case that $m=1$, $H_1^t H_{(2^a1^b)}[X;q,t]$ is also Schur positive,
but not because any of the lemmas in this section hold true and so the following proof
does not apply.
There is an algebraic explanation for this fact.  By using commutation relations
(\ref{comm1}) and (\ref{comm2}) and Theorem \ref{vop2} we have that
$H_1^t H_{(2^a1^b)}[X;q,t] =  t^a H_{(2^a1^{b+1})}[X;qt^{-1},t]$.

\bop
Set $n=2a+b$.  For $m \geq 0$, the action of $H_m^t$ on the Schur function basis is given
in Proposition \ref{snakeruleHm} and using equation (\ref{H2a1btab}), we have that
\begin{equation} \label{HmH2a1b}
H_m^t H_{(2^a1^b)}[X;q,t]=\sum_T \sum_{\rho \slash \la(T) \in \H_{m+n}} 
(-1)^{ht_k(\rho)-1} q^{b_{(2^a1^b)}(T)} 
t^{a_{(2^a1^b)}(T)+| \la \slash \rho^r |} 
s_{\rho \snake_n}[X]
\end{equation}

There is exactly one term in equation (\ref{HmH2a1b}) for every
pair $(T,\rho)$ where $T$ is a standard tableau of size $n$ and $\rho$ is a
partition such that $\rho \slash \la(T) \in \H_{m+n}$ (exactly the pairs for which
$\AHr{m}{\rho} T$ is defined).  The weight of a pair 
$(T,\rho)$ will be $(-1)^{ht_k(\rho)-1} q^{b_{(2^a1^b)}(T)} 
t^{a_{(2^a1^b)}(T)+| \la \slash \rho^r |} 
s_{\rho \snake_n}[X]$.  Note that the immaterial pairs have weight zero.
We will use the lemmas above to show that the
unstable pairs cancel and we will be left with the stable pairs which are in
bijection with the standard tableau containing $\def\Tscale{.65}\tableau[Y]{1,2, ,m}$
as a subtableau.

If $(T,\rho)$ is an unstable pair then let $(\hat{T},\tilde{\rho}) = {\mathcal I}^{m,n}(T,\rho)$.
Since the $(2^a1^b)-type$s are the same then we  have that $b_{(2^a1^b)}(T) = b_{(2^a1^b)}(\hat{T})$
and by Proposition \ref{stats} and Lemma \ref{char},  
\begin{align}
a_{(2^a1^b)}(T) + |\la(T) \slash \rho^r| &=
c(T) + |\la(T) \slash \rho^r|  \nonumber \\
&\hskip .3in - \sum_{i=1}^a ((n+1) - 2i) \chi(type_{(2^a1^b)}(T)_i = \tot) \\
&= c(\AHr{m}{\rho} T ) - \left(\ontop{m}{2}\right)
- (m-1)n   \nonumber \\
&\hskip .3in-\sum_{i=1}^a ((n+1) - 2i) \chi(type_{(2^a1^b)}(T)_i = \tot) \\
&= c(\AHr{m}{\tilde{\rho}} \hat{T} ) - \left(\ontop{m}{2}\right)
- (m-1)n \nonumber\\
&\hskip .3in- \sum_{i=1}^a ((n+1) - 2i) \chi(type_{(2^a1^b)}(\hat{T})_i = \tot) \\
&= c(\hat{T}) + |\la(\hat{T})\slash \tilde{\rho}^r| \nonumber \\
&\hskip .3in-\sum_{i=1}^a ((n+1) - 2i) \chi(type_{(2^a1^b)}(\hat{T})_i = \tot) \\
&= a_{(2^a1^b)}(\hat{T}) + |\la(\hat{T})\slash \tilde{\rho}^r|
\end{align}
Therefore the unstable pairs $(T,\rho)$ and $(\hat{T},\tilde{\rho})$ have 
the same weight but opposite sign.

If $(T,\rho)$ is a stable pair then by Lemma \ref{tfae} the sign of the term
is positive and hence are the only terms
that contribute to this sum.  The map that sends $(T,\rho)$ to $\AHr{m}{\rho}T$ is
a bijection from these pairs to the standard tableaux that contain $\def\Tscale{.65}\tableau[Y]{1,2, ,m}$
as a subtableau.  Also remark that by Proposition \ref{stats} and 
Lemma \ref{char}, for a tableau of this sort
\begin{align}
a_{(2^a1^b)}(\Hti{m} T) + |\la(T) \slash \rho^r| &=
c(\Hti{m} T) + |\la(T) \slash \rho^r|  \nonumber \\
&\hskip .3in - \sum_{i=1}^a ((n+1) - 2i) \chi(type_{(2^a1^b)}(\Hti{m} T)_i = \tot) \\
&= c(T)   - \left(\ontop{m}{2}\right) - (m-1)n    \nonumber \\
&\hskip .3in - \sum_{i=1}^a ((n+1) - 2i) \chi(type_{(2^a1^b)}(\Hti{m} T)_i = \tot) 
\end{align}
Hence, the exponent of $t$ in the weight of the pair
agrees with our definition of $a_{(m2^a1^b)}(T)$.
\eop

\begin{prop} \label{barHmstats}
\begin{equation*}
{\bar H}_m^t H_{(2^a1^b)}[X;q,t] = \sum_T q^{b_{(m2^a1^b)}(T)-\left( \ontop{m}{2} \right)}
t^{a_{(m2^a1^b)}(T)} s_{\la(T)}[X]
\end{equation*}
where the sum is over standard tableau $T$ of size $2a+b+m$ that contain $\def\Tscale{.65}
\tableau[Y]{m| |2|1}$ as
a subtableau and for these tableaux, $ b_{(m2^a1^b)}(T) = b_{(2^a1^b)}(\Hti{m}T)
+\left( \ontop{m}{2} \right)$ and
\begin{equation*}
a_{(m2^a1^b)}(T) = c(T) - \sum_{i=1}^a ((2a+b+1) - 2i) \chi(type_{(2^a1^b)}(\Hti{m}T)_i = \tot)
\end{equation*}
\end{prop}

\bop
All of the arguments and lemmas in this section can be transposed and traced through for the
case of ${\bar H}_m^t$ and $\AVr{m}{\rho}$,
but  an algebraic method is much easier.  Set $n=2a+b$.
We have that
$H_\mu[X;q,t] = \omega t^{n(\mu)} q^{n(\mu')} H_\mu[X;1/q,1/t]$  where
$n(\mu) = \sum_i (i-1) \mu_i$ by \cite[eq. 8.14, p.354]{M}, and hence

\begin{align}
{\bar H}_m^t H_{(2^a1^b)}[X;q,t] &=
\omega H_m^{1/t} \omega R^t H_{(2^a1^b)}[X;q,t] \\
&= \omega H_m^{1/t} t^n t^{n(2^a1^b)} q^{n(a+b,a)} H_{(2^a1^b)}[X;1/q,1/t] \\
&=  t^{n+n(2^a1^b)} q^a \sum_T t^{-a_{(m2^a1^b)}(\omega T)} q^{-b_{(m2^a1^b)}(\omega T)}
s_{\la(T)}[X] \\
&=\sum_T t^{n+n(2^a1^b)- a_{(m2^a1^b)}(\omega T)} q^{b_{(m2^a1^b)}(T)-\left( \ontop{m}{2} \right)}
s_{\la(T)}[X] \label{inky}
\end{align}
where the sum here is over all standard tableaux of size $n+m$ that contain the subtableau 
$\def\Tscale{.65} \tableau[Y]{m| |2|1}$.

The $t$ exponent can be reduced after this using Proposition \ref{Hmstats},
Fact \ref{fact3}, and the relation $a_{(2^a1^b)}(\omega T) =
n(2^a1^b) - a_{(2^a1^b)}(T)$.
\begin{align}
n+n(2^a1^b)- a_{(m2^a1^b)}(\omega T) &=
n(2^a1^b)- c(\omega T) + \left( \ontop{m}{2} \right) + m n  \nonumber\\
&\hskip .3in+ \sum_{i=1}^a ((n+1) - 2i) \chi(type_{(2^a1^b)}(\omega \Hti{m}T)_i = \tot) \\
&= c(T) - \left( \ontop{n}{2} \right) + n(2^a1^b) \nonumber\\
&\hskip .3in+ \sum_{i=1}^a ((n+1) - 2i) \chi(type_{(2^a1^b)}(\omega \Hti{m}T)_i = \tot) \\
&= c(T) - \sum_{i=1}^a ((n+1) - 2i) \chi(type_{(2^a1^b)}(\Hti{m}T)_i = \tot)
\end{align}
\eop

To extend these results to $\Hbac$ and $\Hcab$, recovering the statistics involves deleting
the first cell and then reducing to the two column case.  We introduce the notation of the
operator $K_h$ that acts on standard tableaux with more than $h$ cells.  The result $K_h T$ is
the tableau formed by deleting the cells with labels $1$ through $h$ and then lowering the 
labels in the remaining tableau by $h$.

\begin{prop} \label{Hbacstats}
\begin{equation*}
\Hbac H_{(2^a1^b)}[X;q,t] = \sum_T q^{b_{(32^a1^b)}(T)-1}
t^{a_{(32^a1^b)}(T)} s_{\la(T)}[X]
\end{equation*}
where the sum is over standard tableau $T$ of size $2a+b+3$ that contain $\tbac$ as
a subtableau. For these tableaux set 
$b_{(32^a1^b)}(T) = b_{(2^a1^b)}(\Hti{2}K_1 T)+1$ and
\begin{equation*}
a_{(32^a1^b)}(T) = c(T) - 1 -(2a+b) - \sum_{i=1}^a ((2a+b+1) - 2i) \chi(type_{(2^a1^b)}(\Hti{2}K_1 T)_i = \tot)
\end{equation*}
\end{prop}

\bop
\begin{equation}
e_1[X] H_2^t H_{(2^a1^b)}[X;q,t] = \sum_T q^{b_{(2^{a+1}1^b)}(K_1 T)}
t^{a_{(2^{a+1}1^b)}(K_1 T)} s_{\la(T)}[X]
\end{equation}
where the sum is over all standard tableaux $T$ that contain either $\tabc$ or $\tbac$.

For the tableaux that contain $\tabc$, $c(T) = c(K_1 T) + 2a+b+2$ and 
\begin{align}
a_{(2^{a+1}1^b)}(K_1 T) 
&= c(K_1 T) - \sum_{i=1}^{a+1} ((2a+b+3) - 2i) \chi(type_{(2^{a+1}1^b)}(K_1 T)_i = \tot) \\
&= c(T) -(2a+b+2)-(2a+b+1) \nonumber\\
&\hskip .3in- \sum_{i=1}^{a+1} ((2a+b+1) - 2i) 
\chi(type_{(2^{a}1^b)}(\Hti{2} K_1 T)_i = \tot) \label{gooby}
\end{align}

For those tableau that have $\Hti{3} T = \Hti{2} K_1 T$, it is
clear that equation (\ref{gooby}) is equal to $a_{(32^a1^b)}(T)$.  For the
rest we need verify that $type_{(2^a1^b)}(\Hti{2}K_1 T) = type_{(2^a1^b)}( \Hti{3} T)$.

Let $\bar{T} = \Hti{3} T$ so that $T = \AHr{3}{(3+2a+b+\la(T)_1,\la(T)^r)} \bar{T}$.
Therefore we have that 
$K_1 T = K_1 \AHr{3}{(3+2a+b+\la(T)_1,\la(T)^r)} \bar{T} = \AHr{2}{(2+2a+b+\la(T)_1,\la(T)^r)}
\bar{T}$.
\begin{equation} \label{tpeq}
type_{(2^a1^b)}(\Hti{2}K_1 T) = type_{(2^a1^b)}(\Hti{2}\AHr{2}{(2+2a+b+\la(T)_1,\la(T)^r)} \bar{T}) =
type_{(2^a1^b)}(\bar{T})
\end{equation}

From this we can conclude that
\begin{equation}
(e_1[X] H_2^t - H_3^t) H_{(2^a1^b)}[X;q,t] = \sum_T q^{b_{(2^{a+1}1^b)}(K_1 T)}
t^{a_{(2^{a+1}1^b)}(K_1 T)} s_{\la(T)}[X]
\end{equation}
where the sum is over all standard tableaux $T$ that contain $\tbac$.
For these tableaux we have that $c(T) = c(K_1 T)$ and hence
\begin{align}
a_{(2^{a+1}1^b)}(K_1 T) 
&= c(K_1 T) - \sum_{i=1}^{a+1} ((2a+b+3) - 2i) \chi(type_{(2^{a+1}1^b)}(K_1 T)_i = \tot) \\
&= c(T) - (2a+b+1) - \sum_{i=1}^{a} ((2a+b+1) - 2i) \chi(type_{(2^{a}1^b)}(\Hti{2} K_1 T)_i = \tot)
\end{align}

\eop

A similar proof verifies that the operator $\Hcab$ works as expected.

\begin{prop} \label{Hcabstats}
\begin{equation*}
\Hcab H_{(2^a1^b)}[X;q,t] = \sum_T q^{b_{(32^a1^b)}(T)-2}
t^{a_{(32^a1^b)}(T)} s_{\la(T)}[X]
\end{equation*}
where the sum is over standard tableau $T$ of size $2a+b+3$ that contain $\tcab$ as
a subtableau. For these tableaux set 
$b_{(32^a1^b)}(T) = b_{(2^a1^b)}(\Hti{2}\bar{T})+2$ and
\begin{equation*}
a_{(32^a1^b)}(T) = c(T) - 2 -(2a+b) - \sum_{i=1}^a ((2a+b+1) - 2i) \chi(type_{(2^a1^b)}(\Hti{2}K_1 T)_i = \tot)
\end{equation*}
\end{prop}

These propositions together as well as Theorem \ref{vop3} and equation (\ref{nice3}) show
the following theorem.

\begin{thm}
\begin{equation*}
H_{(32^a1^b)}[X;q,t] = \sum_T q^{b_{(32^a1^b)}(T)}
t^{a_{(32^a1^b)}(T)} s_{\la(T)}[X]
\end{equation*}
where the sum is over all standard tableaux $T$ of size $3+2a+b$ and $a_{(32^a1^b)}(T)$
and $b_{(32^a1^b)}(T)$ are defined in Propositions
\ref{Hmstats}, \ref{barHmstats}, \ref{Hbacstats} and \ref{Hcabstats}.
\end{thm}

Unfortunately, since there is no general formula for the $H^S$
(each one is a special case) and each must be checked separately
to show that they have the property that $H^S H_{(2^a1^b)}[X;q,t]$ is
Schur positive.

For $T$ a standard tableau of size $2a+b+4$ that contains the standard tableau $S$ of
size $4$ we will define that

\begin{equation} \label{a42a1b}
a_{(42^a1^b)}(T) = c(T) - \alpha_S - \beta_S (2a+b) 
- \sum_{i=1}^a ((2a+b+1) - 2i) \chi(type_{(2^a1^b)}(\theta_S(T))_i= \tot)
\end{equation}

\begin{equation} \label{b42a1b}
b_{(42^a1^b)}(T) = b_{(2^a1^b)}(\theta_S(T))+\gamma_S
\end{equation}
where $\alpha_S$, $\beta_S$, $\theta_S(T)$, and $\gamma_S$ are given in Table 1.

The format for the statistics is the same, just a few variables change for
each equation.  It is easier to give a table of values for the parts of
this equation that change than it is to describe a procedure for calculating
them. Below is a table for the statistics for $H_{(42^a1^b)}[X;q,t]$
(The values in the table for $H_{(32^a1^b)}[X;q,t]$ follow from Propositions
\ref{Hmstats}, \ref{barHmstats}, \ref{Hbacstats} and \ref{Hcabstats}).

\begin{table}[tbp]
\begin{tabular}[b]{cccccl}
S & $\alpha_S$ & $\beta_S$ & $\theta_S(T)$ & $\gamma_S$ & $H^S$\\ \hline \\
$\tabcd$ & $6$ & $3$ & $\Hti{4} T$     & $0$ & $H_4^t$ \\ \\
$\tbacd$ & $3$ & $2$ & $\Hti{3} K_1 T$ & $1$ & $e_1 H_3^t - H_4^t$ \\ \\
$\tcabd$ & $4$ & $2$ & $\Hti{2} K_2 T$ & $2$ & $h_2 H_2^t - H_4^t$  \\ \\
$\tdabc$ & $5$ & $2$ & $\Hti{2} K_2 T$ & $3$ & $h_2 {\bar H}_2^t - h_1 {\bar H}_3^t +  {\bar H}_4^t$ \\ \\
$\tcdab$ & $4$ & $2$ & $\Hti{2} K_2 T$ & $2$ & $h_2 H_2^t - H_4^t$ \\ \\
$\tbdac$ & $2$ & $1$ & $\Hti{2} K_2 T$ & $4$ & $e_2 {\bar H}_2^t - {\bar H}_4^t$ \\ \\
$\tcbad$ & $1$ & $1$ & $\Hti{2} K_2 T$ & $3$ & $e_2 H_2^t - e_1 H_3^t +  H_4^t$ \\ \\
$\tdbac$ & $2$ & $1$ & $\Hti{2} K_2 T$ & $4$ & $e_2 {\bar H}_2^t - {\bar H}_4^t$ \\ \\
$\tdcab$ & $3$ & $1$ & $\Hti{3} K_1 T$ & $5$ & $e_1 {\bar H}_3^t - {\bar H}_4^t$ \\ \\
$\tdcba$ & $0$ & $0$ & $\Hti{4} T$     & $6$ & ${\bar H}_4^t$ \\ \\ \hline \\
$\tabc$ & $3$ & $2$ & $\Hti{3} T$      & $0$ & $H_3^t$ \\ \\
$\tbac$ & $2$ & $1$ & $\Hti{2} K_1 T$  & $1$ & $e_1 H_2^t - H_3^t$ \\ \\
$\tcab$ & $1$ & $1$ & $\Hti{2} K_1 T$  & $2$ & $e_1 \bar{H}_2^t - \bar{H}_3^t$ \\ \\
$\tcba$ & $0$ & $0$ & $\Hti{3} T$      & $3$ & $\bar{H}_3^t$ \\ \\ \hline
\end{tabular}
\caption{Values for variable pieces of formulas (\ref{a42a1b}) and (\ref{b42a1b})}
\end{table}

\begin{prop} \label{H4stats}
For each standard tableau $S$ of size $4$ except $\tcabd$, $\tcdab$, $\tbdac$ and
$\tdbac$ we have that
\begin{equation*}
H^S H_{(2^a1^b)}[X;q,t] = \sum_T q^{b_{(42^a1^b)}(T) - \gamma_S} t^{a_{(42^a1^b)}(T)} s_{\la(T)}[X]
\end{equation*}
where the sum is over all standard tableaux $T$ of size $4+2a+b$ that contain $S$ as a
subtableau.  In addition we have that
\begin{equation*}
\left(\Hcdab+\Hcabd \right) H_{(2^a1^b)}[X;q,t] = \sum_T q^{b_{(42^a1^b)}(T) - 2} t^{a_{(42^a1^b)}(T)} s_{\la(T)}[X]
\end{equation*}
where the sum is over all standard tableaux $T$ of size $4+2a+b$ that contain either
$\tcabd$ or $\tcdab$.  Furthermore,
\begin{equation*}
\left(\Hbdac+\Hdbac\right) H_{(2^a1^b)}[X;q,t] = \sum_T q^{b_{(42^a1^b)}(T) - 4} t^{a_{(42^a1^b)}(T)} s_{\la(T)}[X]
\end{equation*}
where the sum is over all standard tableaux $T$ of size $4+2a+b$ that contain either
$\tbdac$ or $\tdbac$.
\end{prop}

The fact that $H_{(42^a1^b)}[X;q,t]$ is Schur positive follows immediately from this proposition
and Theorem \ref{vop4} and equation (\ref{expv4}).

\begin{thm}
\begin{equation*}
H_{(42^a1^b)}[X;q,t] = \sum_T q^{b_{(42^a1^b)}(T)}
t^{a_{(42^a1^b)}(T)} s_{\la(T)}[X]
\end{equation*}
where the sum is over all standard tableaux $T$ of size $4+2a+b$.
\end{thm}

\bop (of Proposition \ref{H4stats})

We note that if $T$ contains $\tabcd$ or $\tdcba$, it is not necessarily the case that
 $\Hti{4} T = \Hti{3} K_1 T = \Hti{2} K_2 T$, but it will be true that 
$type_{(2^a1^b)}(\Hti{4} T) = type_{(2^a1^b)}(\Hti{3} K_1 T) =type_{(2^a1^b)}(\Hti{2} K_2 T)$.
Similarly, if $T$ contains $\tdcab$ or $\tbacd$ then in general it is not true
that $\Hti{3} K_1 T = \Hti{2} K_2 T$, but it will be true that
$type_{(2^a1^b)}(\Hti{3} K_1 T) =type_{(2^a1^b)}(\Hti{2} K_2 T)$.  This follows from a
similar argument to the one used in equation (\ref{tpeq}).

For each operator $H^S$ we must verify that it has the property stated
in the propostion.  The case of $\Habcd$ and $\Hdcba$
follow from Proposition \ref{Hmstats} and \ref{barHmstats} respectively.

\vskip .2in
\noindent
$\Hbacd$:

Since $H_3^t H_{(2^a1^b)}[X;q,t]$ is a generating function for the standard tableaux
that contain $\tabc$ then $e_1[X] H_3^t H_{(2^a1^b)}[X;q,t]$ is a generating function
that for the standard tableaux that contain $\tabcd$ or $\tbacd$, that is
\begin{equation}
e_1[X] H_3^t H_{(2^a1^b)}[X;q,t] = \sum_T q^{b_{(2^a1^b)}(\Hti{3} K_1 T)}
t^{a_{(32^a1^b)}(K_1 T)} s_{\la(T)}[X]
\end{equation}
To verify this
proposition in this case we need to show that
$a_{(32^a1^b)}(K_1 T) = a_{(42^a1^b)}(T)$ for these tableaux.  For the tableaux that contain
$\tbacd$, this is true by definition.  For the tableaux that contain $\tabcd$ we note
that $c(K_1 T) = c(T) - (2a+b+3)$ and $type_{(2^a1^b)}(\Hti{4} T) = type_{(2^a1^b)}(\Hti{3} K_1 T)$.

\vskip .2in
\noindent
$\Hcdab + \Hcabd$:

$h_2 H_2^t H_{(2^a1^b)}[X;q,t]$ is a generating function for the standard tableaux that contain
$\tcdab$, $\tcabd$ or $\tabcd$ as subtableau.
\begin{equation}
h_2[X] H_2^t H_{(2^a1^b)}[X;q,t] = \sum_T q^{b_{(2^a1^b)}(\Hti{2} K_2 T)}
t^{a_{(2^{a+1}1^b)}(K_2 T)} s_{\la(T)}[X]
\end{equation}
We verify for these standard tableaux that
$a_{(2^{a+1}1^b)}(K_2 T) = a_{(42^a1^b)}(T)$
For the tableaux that contain $\tcdab$ or $\tcabd$ as standard subtableau we have that
$c(K_2 T) = c(T) - (2a+b+3)$.  For the tableaux that contain $\tabcd$, $c(K_2 T) = c(T) -
(2(2a+b)+5)$ and we note that $type_{(2^a1^b)}(\Hti{3} K_1 T) = type_{(2^a1^b)}(\Hti{2} K_2 T)$.

\vskip .2in
\noindent
$\Hcbad$:

$e_2 H_2^t H_{(2^a1^b)}[X;q,t]$ is a generating function for the standard tableaux that
contain $\tcbad$ and $\tbacd$ as a subtableau with the formula
\begin{equation}
h_2[X] H_2^t H_{(2^a1^b)}[X;q,t] = \sum_T q^{b_{(2^a1^b)}(\Hti{2} K_2 T)}
t^{a_{(2^{a+1}1^b)}(K_2 T)} s_{\la(T)}[X]
\end{equation}
The operator $\Hcbad = e_2 H_2^t - \Hbacd$.  We calculate that
for the standard tableaux that contain
$\tcbad$, $c(K_2 T) = c(T)$ and for the standard tableaux the contain
$\tbacd$, $c(K_2 T) = c(T) - (2a+b+2)$.

\vskip .2in
The verification of the cases for the operators with the transpose tableaux is nearly
identical or an algebraic argument can be used.  $\alpha_S = 6 - \alpha_{\omega S}$
and $\beta_S = 3 -\beta_{\omega S}$ for all standard tableaux $S$ of size $4$.
Let $n=2a+b$.
It follows from the relation
$H^S = \omega H^{\omega S} \coeff_{t\rightarrow 1/t} \omega R^t$ and
equation (\ref{inky}) that 

\begin{equation}
H^{\omega S}H_{(2^a1^b)}[X;q,t] = \sum_T t^{n+n(2^a1^b)- a_{(42^a1^b)}(\omega T)}
q^{n(a+b,a)-b_{(42^a1^b)}(\omega T)} s_{\la(T)}[X]
\end{equation}
where the sum is over all standard tableaux $T$ that contain $\omega S$.

Now for the statistic, we have that

\begin{align}
n+n(2^a1^b)- a_{(42^a1^b)}(\omega T) &= n+n(2^a1^b)- (c(\omega T) - \alpha_{S} - \beta_{S} n \\
&\hskip .3in- \sum_{i=1}^{a} ((2a+b+1) - 2i) \chi(type_{(2^{a}1^b)}(\theta_{S}(\omega T))_i = \tot)) \nonumber \\
&= c(T) -(n+4)(n+3)/2+n(n-1)/2+n + \alpha_{S} + \beta_{S} n \\
&\hskip .3in - \sum_{i=1}^{a} ((2a+b+1) - 2i) \chi(type_{(2^{a}1^b)}(\theta_{\omega S}(T))_i = \tot) \nonumber \\
&= c(T) - (3-\beta_S)n - (6-\alpha_S) \\
&\hskip .3in - \sum_{i=1}^{a} ((2a+b+1) - 2i) \chi(type_{(2^{a}1^b)}(\theta_{\omega S}(T))_i = \tot) \nonumber \\
&= a_{(42^a1^b)}(T) 
\end{align}
\eop

It makes sense here to generalize the notion of the $type_{(2^a1^b)}$ as it was introduced
in \cite{Za1} to the case of $(32^a1^b)$ and $(42^a1^b)$.  For $m = 3$ or $4$, $T$ be a 
standard tableau of size $m+2a+b$ that contains the standard tableau $S$ of size $m$.
Define the $type_{(m2^a1^b)}(T) = (S, type_{(2^a1^b)}(\theta_S(T)))$ where $\theta_S(T)$
is given in Table 1.

We finish this article by noting that the following analog to Conjecture 4.8 of \cite{Za1}
seems to be true by experimental calculation.

\begin{conjecture}
The number of tableaux of a fixed $(m2^a1^b)-type$, $(S,s)$, 
and fixed $a_{(m2^a1^b)}$ value, $i$, increases to
a maximum and then decreases as $i$ ranges over all possible values.

More precisely stated, if $S$ is a standard tableau of size $m$ ($m$ = $3$ or $4$) and
$s$ is a sequence of $a$ standard tableaux of size $2$ and $b$ of size $1$,
then let $A_{(S,s)}^i = \# \{ T | T \in ST^{m+2a+b}, type_{(m2^a1^b)}(T) = (S,s),
a_{(m2^a1^b)}(T) = i\}$.  The sequence $A_{(S,s)}^* = \left( A_{(S,s)}^0,  A_{(S,s)}^1,
A_{(S,s)}^2, \ldots \right)$ is unimodal.
\end{conjecture}

\begin{example}
Let $s = (\to, \to,\to)$ and let $S$ range over all standard tableaux of size $3$.
\begin{align*}
S=\tabc \hskip .4in & A_{(S,s)}^* = (1,2,3,4,2,1,1)\\
S=\tbac \hskip .4in & A_{(S,s)}^* = (2,4,6,5,4,2,1)\\
S=\tcab \hskip .4in & A_{(S,s)}^* = (1,2,4,5,6,4,2)\\
S=\tcba \hskip .4in & A_{(S,s)}^* = (1,1,2,4,3,2,1)
\end{align*}

Let $s=(\to,\to)$ and let $S$ range over all standard tableaux of size $4$.  This
partitions all of the standard tableaux of size $6$ into 10 different types.
\begin{equation*}
\begin{array}{llll}
S=\tabcd \hskip .4in & A_{(S,s)}^* = (1,2,1,1)\hskip .5in& 
S=\tbacd \hskip .4in & A_{(S,s)}^* = (2,4,2,1)\\\\
S=\tcabd \hskip .4in & A_{(S,s)}^* = (2,4,2,1)&
S=\tdabc \hskip .4in & A_{(S,s)}^* = (1,2,3,3)\\\\
S=\tcdab \hskip .4in & A_{(S,s)}^* = (1,2,2,1)&
S=\tbdac \hskip .4in & A_{(S,s)}^* = (1,2,2,1)\\\\
S=\tcbad \hskip .4in & A_{(S,s)}^* = (3,3,2,1)&
S=\tdbac \hskip .4in & A_{(S,s)}^* = (1,2,4,2)\\\\
S=\tdcab \hskip .4in & A_{(S,s)}^* = (1,2,4,2)&
S=\tdcba \hskip .4in & A_{(S,s)}^* = (1,1,2,1)
\end{array}
\end{equation*}
\end{example}

\begin{example}
Let $s = (\tot,\tot,\to)$ and $S$ range over all standard tableaux of size $3$.

\begin{align*}
S = \tabc \hskip .4in & A_{(S,s)}^* = (1,4,6, 8, 7, 6, 4,2,1,1) \\
S = \tbac \hskip .4in & A_{(S,s)}^* = (2,3,8,12,13,10, 8,4,2,1) \\
S = \tcab \hskip .4in & A_{(S,s)}^* = (0,0,2, 7,12,14,13,9,4,2)\\
S = \tcba \hskip .4in & A_{(S,s)}^* = (0,0,1, 3, 5, 7, 7,4,2,1)
\end{align*}
\end{example}

\end{document}